\RequirePackage{fix-cm}
\documentclass[smallcondensed]{svjour3}     
\smartqed  
\usepackage{graphicx}
\usepackage{amssymb}
\usepackage{hhline}
\usepackage{amsmath}
\usepackage{graphicx}
\usepackage{amsfonts}
\usepackage{color}
\usepackage{array}
\usepackage{subfig}

\usepackage{pgfplots}
\usepackage{pgfplotstable}
\definecolor{markercolor}{RGB}{124.9, 255, 160.65}
\pgfplotsset{width=10cm,compat=1.3}
\pgfplotsset{
tick label style={font=\small},
label style={font=\small},
legend style={font=\small}
}

\usepackage{filecontents}
\usepackage{verbatim}

\usepackage{pgfplots}
\usepackage{pgfplotstable}
\usetikzlibrary{patterns}
\definecolor{markercolor}{RGB}{124.9, 255, 160.65}
\pgfplotsset{width=10cm,compat=1.3}
\pgfplotsset{
tick label style={font=\scriptsize	},
label style={font=\scriptsize	},
legend style={font=\scriptsize}
}

\newcommand{\bigo}{\color{black}\mathcal{O}\color{black}}

\newcommand{\mb}[1]{\mathbf{#1}}

\newcommand{\LRp}[1]{\left( #1 \right)}
\newcommand{\LRs}[1]{\left[ #1 \right]}

\begin{document}

\title{Leapfrog time-stepping for Hermite methods
}


\author{Arturo Vargas         \and
        Thomas Hagstrom \and \\
         Jesse Chan \and
        Tim Warburton 
}


\institute{A. Vargas \at
             Lawrence Livermore National Laboratory,
             \\ 7000  East Ave, Livermore, CA 94550 \\
              \email{vargas45@llnl.gov}           
           \and
           T. Hagstrom \at
           Department of Mathematics, Southern Methodist University, Dallas, TX 75275
           \and
           J. Chan \at  Department of Computational and Applied Mathematics, Rice University, 6100 Main Street, Houston, TX 77005
           \and 
           T. Warburton \at Department of Mathematics, Virginia Tech, Blacksburg, VA 24060
}

\date{Received: date / Accepted: date}

\maketitle

\begin{abstract}
We introduce Hermite-leapfrog methods for first order wave systems. The new Hermite-leapfrog methods 
pair leapfrog time-stepping with the Hermite methods of Goodrich and co-authors \cite{goodrich2006hermite}. The new schemes stagger field variables in both time and space and are high-order accurate. We provide a detailed description of the method and demonstrate that the method conserves variable quantities in one-space dimension. Higher dimensional versions of the method are constructed via a tensor product construction. Numerical evidence and rigorous analysis in one space dimension establish stability and high-order convergence. Experiments demonstrating efficient implementations on a graphics processing unit are also presented. 
\end{abstract}

\section{Introduction}
Simulations of wave propagation play a crucial role in science and engineering. In
applications to geophysics, they are the engine of many seismic imaging algorithms.
For electrical engineers, they can be a useful tool for the design of radars and antennas. In these applications achieving high fidelity simulations is challenging due to the inherent issues in modeling highly oscillatory waves and the associated high computational cost of high-resolution simulations. Thus the ideal numerical method should be able to capture high-frequency waves and be suitable for parallel computing.

The staggered finite-difference time-domain (FDTD) method has been widely adopted as the method of choice for simulating electromagnetic and acoustic wave propagation on structured grids \cite{bencomo2015discontinuous,kowalczyk2010comparison}. In the context of electromagnetics the method was first pioneered by K.S. Yee in \cite{yee1966numerical}; while the discretization for the acoustic wave equation was introduced by Virieux in \cite{gauthier1986two} and Levander in \cite{levander1988fourth}. The underlying idea of staggered FDTD methods is to stagger field variables both in time and space and propagate the solution using leapfrog time-stepping \cite{leveque2007finite}. This type of discretization is known to be efficient and easy to implement but is typically only second order accurate and struggles to resolve highly-oscillatory waves. Fourth order accurate variants of the method have been proposed in the literature \cite{xie2002explicit,yefet2001staggered} but our work differs as we present methods of arbitrary order.

In contrast to FDTD methods where derivatives are approximated via difference formulas, the Hermite methods of Goodrich and co-authors \cite{goodrich2006hermite} combine Hermite-Birkhoff interpolation and a staggered (dual) grid to produce high-order numerical methods for first order hyperbolic problems. The resulting ``Hermite methods" carry out a localized polynomial reconstruction within cells of a structured grid and propagate the solution completely independent of neighboring cells via Hermite-Taylor time-stepping. The remarkable feature of Hermite methods is that the CFL condition of the resulting scheme is independent of the approximation order \cite{goodrich2006hermite}. A formal introduction to these methods may be found in \cite{goodrich2006hermite,vargas2017hermite,vargas2017variations}. 

Though Hermite methods are relatively new, several variations have been proposed to enhance their applicability and efficiency. To enable geometric flexibility, Hermite methods have been combined with discontinuous Galerkin methods \cite{chen2014hybrid}. Adaptive variants of Hermite methods have been introduced by means of $p$-adaptivity in \cite{chen2012p}, and preliminary work on $h$-adaptivity has been presented in \cite{appelo2011hermite}. Variations of the method which do not require a dual grid have been introduced in \cite{vargas2017variations}, and a flux conservative Hermite method has been introduced in \cite{kornelus2018flux}.

Inspired by the favorable features of Hermite methods, this work presents a new variant which combines leapfrog time stepping and staggered solutions in space and time. The resulting ``Hermite-leapfrog" scheme may be viewed as a high order variant of the Yee scheme. We present this work as complimentary to the recent work of Appel\"{o} and co-authors \cite{appelo2018hermite} in which a similar time-stepping scheme was introduced for the second order acoustic wave equation. In this work we focus on first order wave systems and note that an early version appeared in the thesis of Vargas \cite{vargas2017hermite}.

The remainder of this article proceeds as follows. We start with a detailed description in one space dimension. We then demonstrate that the method conserves variable quantities, establishing stability, and prove general error estimates in one space dimension. In addition we look at the dispersion relation, which provides an explanation for aspects of the method's accuracy observed in the numerical experiments. Experimental convergence rates are reported for both constant and spatially varying wave speed.  In two-dimensions we discuss incorporating zero Dirichlet boundary conditions and demonstrate the method's ability to resolve highly oscillatory waves. Lastly, we demonstrate computational efficiency on a graphics processing unit.

\section{Description of the method}
To introduce the method, we consider the following one-dimensional wave system
\begin{align} \label{eq:waveModel}
\frac{\partial p}{\partial t} &=  c \frac{\partial v}{\partial x},  \quad\frac{\partial v}{\partial t} = c \frac{\partial p}{\partial x} \\ \nonumber
 p(x,t_0) &= f(x), \quad v(x,t_{0}+{\Delta t}/{2}) = g(x).  \nonumber
\nonumber
\end{align}
In this example, $c$ corresponds to the speed of the propagating wave. The degrees of freedom of the Hermite-leapfrog method are the function value and first $m$ derivatives at the nodes of a structured grid. An $m^{th}$ order scheme discretizes the equations by staggering pressure and velocity approximations in both space and time. The discretization of the pressure occurs on a primary grid $\Omega$ while the discretization of the velocity term is maintained on a staggered dual grid $\tilde{\Omega}$. For the purpose of introducing the method we assume that initial conditions are given at grid-points staggered in both time and space. We define the primary grid as a collection of $K$ equidistant points

\begin{align}
\Omega = \{x_j: \quad x_j = x_{\min} + j h_x, \quad j = 0,\ldots, K-1 \} \label{eq:pGrid},
\end{align}
analogously the dual grid, which holds velocity approximations, is defined as
\begin{align}
\tilde{\Omega} = \{ x_{j+1/2} = x_{\min} + \left(j + \frac{1}{2} \right) h_x, \quad j = 0,\ldots, K-1 \} \label{eq:dGrid} .
\end{align}
Approximations of each variable are carried out through local Taylor series expansions at each node $x_j$. The expansion takes the form 
\[
u(x) \approx u_j(x) = \sum_{i=0}^{m} \tilde{u}_i \left( \frac{ x-x_j}{h} \right)^i,
\]
where $\tilde{u}_i$ corresponds to scaled approximations of the $i^{th}$ derivative of $u(x)$
\[
\tilde{u}_i = \frac{h^i}{i!}\frac{d^i u}{d x^i}, \quad i=0,\dots m. 
\]

\begin{figure}[h!]
\center
\includegraphics[scale=0.4]{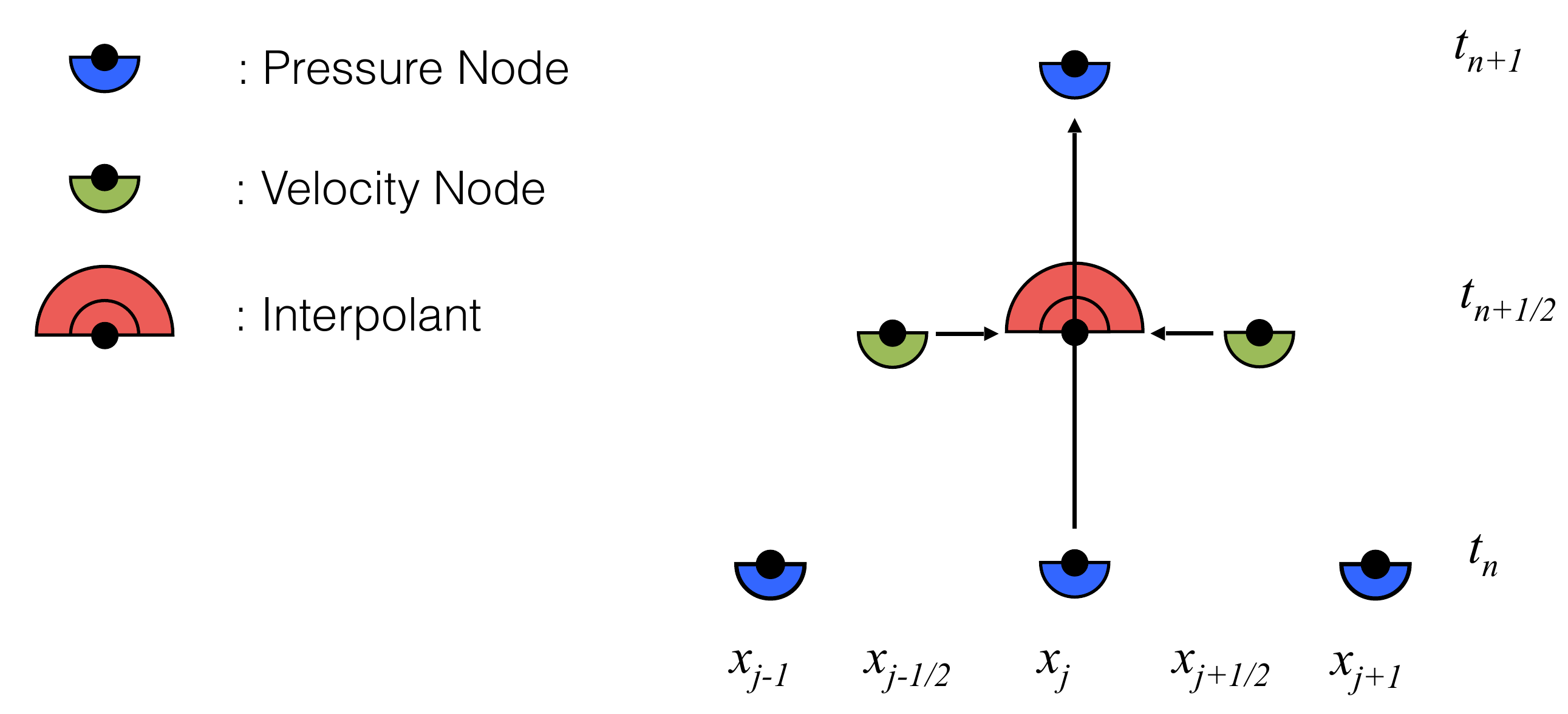}
\caption{The first step of the Hermite-leapfrog method is the use of Hermite-Birkhoff interpolation to approximate the function value and first $2m+1$ derivatives of the velocity field at node $x_j$. The solution of the pressure is then propagated via the Hermite-leapfrog time-stepping algorithm (Equation \ref{eq:final1} ). }\label{fig:plusHermite}
\end{figure}

The time-stepping algorithm of the Hermite-leapfrog scheme is derived by considering the following temporal series at each grid point 
\begin{subequations}
\begin{align}
     p(x,t + \Delta t) &= \sum_{j=0}  \frac{({\Delta t}/{2})^j}{j!} \frac{\partial^{j} p(x, t + \frac{\Delta t}{2}) } {\partial t^{j}}, \label{eq:option1}  \\ 
     p(x,t) &= \sum_{j=0}  \frac{({-\Delta t}/{2})^j}{j!} \frac{\partial^{j} p(x, t + \frac{\Delta t}{2}) }{\partial t^{j}} \label{eq:option2}.
\end{align}
\end{subequations}
Subtracting equations \ref{eq:option1} and \ref{eq:option2} yields
\begin{align}
p(x,t+\Delta t) - p(x,t)  &=  \sum_{j=1,odd}  2 \frac{({\Delta t}/{2})^j}{j!} \frac{\partial^{j} p(x, t+\frac{\Delta t}{2}) }{\partial t^{j}}, \label{eq:update1}
\end{align}
where the limit of the sum is determined by the order of the method. Temporal derivatives are then exchanged for spatial derivatives by means of the Cauchy-Kowalevsky recurrence relation \cite{goodrich2006hermite}; in the case of constant coefficients the recurrence relation simplifies to 
\begin{equation}
\frac{\partial^r p}{\partial t^r} =   c^{r} \frac{\partial^r v}{\partial x^r} , \quad r \in\{1,3,5,\dots\}. \label{eq:LFloop}
\end{equation}
yielding the time-stepping algorithm  
\begin{align} 
\frac{\partial^s p(x,t + \Delta t)}{\partial x^s} &= \frac{\partial^s p (x,t)}{\partial x^s} + 2 \sum_{j=1,odd}  \frac{({c \Delta t}/{2})^j}{j!} \frac{\partial^{j+s} v (x,t + {\Delta t}/{2})}{\partial x^{j+s}}.  \label{eq:final1}
\end{align}
Here we introduce the index $s$ to correspond to the degree of freedom, $s\in \{0, \dots, m\}$, at each grid point.
The key ingredient here is the approximation of the spatial derivatives (the right hand side of Equation \ref{eq:final1}). The approximation is carried out by means of Hermite-Birkhoff interpolation as done in the classic Hermite method \cite{goodrich2006hermite}. The interpolation procedure constructs a polynomial $\tilde{v}_j$  by interpolating the function value and derivatives at $R_j = x_{j+{1}/{2}}$ $L_j = x_{j-{1}/{2}}$
\[
\left.\frac{\partial^i \tilde{v}_j}{\partial x^i}\right|_{L_j} = \left.\frac{\partial^i v_j^L}{\partial x^i}\right|_{L_j}, \qquad 
\left.\frac{\partial^i \tilde{v}_j}{\partial x^i}\right|_{R_j}= \left.\frac{\partial^i v_j^R}{\partial x^i}\right|_{R_j}, \qquad i = 0,\ldots, 2m+1.
\]
Thus, the resulting polynomial reproduces the function value and first $m$ derivatives at the left and right end points. Notably, the coefficients of the polynomial are the approximation of the function value and first $2m+1$ derivatives at the node $x_j$. As demonstrated in 
\cite{vargas2017variations} this results in the following system 
\begin{equation}
\LRs{\begin{array}{c}
\mb{C}^L \\
\mb{C}^R
\end{array}}\tilde{\mb{v}} = \LRs{\begin{array}{c}
\mb{v}^L \\
\mb{v}^R
\end{array}}, 
\label{eq:interp}
\end{equation}
where $\mb{C}^L, \mb{C}^{R}$ is defined as
\begin{align*}
\mb{C}^L_{ls} &= 
\begin{cases}
\LRp{\frac{L_j-\tilde{x}_j}{h}}^{s-l}\frac{1}{l!}\prod\limits_{t=0}^{l-1}(s-t), & s\geq l\\
0, &s < l
\end{cases}\\
\mb{C}^R_{ls} &= 
\begin{cases}
\LRp{\frac{R_j-\tilde{x}_j}{h}}^{s-l}\frac{1}{l!}\prod\limits_{t=0}^{l-1}(s-t), & s\geq l\\
0, &s < l.
\end{cases}
\end{align*}
A full description of the interpolation procedure can be found in  \cite{goodrich2006hermite,vargas2017hermite,vargas2017variations}. The coefficients of the polynomial are used to carry out the time stepping algorithm. Figure \ref{fig:plusHermite} provides an illustration of the stencil associated with the scheme. Next, we describe the stability and convergence properties of the method. 

\section{Stability and convergence}
The essential properties of Hermite-Birkhoff interpolation used to establish stability and to estimate the error are proven in detail in the original paper \cite{goodrich2006hermite}. First and foremost is the orthogonality lemma, which shows that the interpolation process is contracting in a certain Sobolev seminorm. Precisely, let $\mathcal{I}_m$ denote the Hermite interpolation operator employing derivative data up to order $m$; here we will not distinguish between interpolation on the two grids. We will also assume that all functions are periodic and sufficiently smooth. In one space dimension we use the Sobolev semi-inner-product
\begin{equation*}
    \langle f,g \rangle_{m+1} = \int \frac{\partial^{m+1} f}{\partial x^{m+1}} \cdot \frac{\partial^{m+1} g}{\partial x^{m+1}} ,
\end{equation*}
which we generalize in $d$ dimensions to
\begin{equation*}
    \langle f,g \rangle_{m+1} = \int \frac{\partial^{(m+1)d} f}{\partial x_1^{m+1} \cdots \partial x_d^{m+1}} \cdot \frac{\partial^{(m+1)d} g}{\partial x_1^{m+1} \cdots \partial x_d^{m+1}} . 
\end{equation*}
Then, as shown in \cite{goodrich2006hermite}, for any $f$ and $g$ we have
\begin{equation}
    \langle \mathcal{I}_m f , g-\mathcal{I}_m g \rangle_{m+1} =0. \label{ortheq}
\end{equation}
Denoting by $\arrowvert \cdot \arrowvert_{m+1}$ the seminorm induced by these semi-inner-products, we thus have, by the Pythagorean Theorem
\begin{equation}
    \arrowvert f \arrowvert_{m+1}^2 = \arrowvert \mathcal{I}_m f \arrowvert_{m+1}^2 + \arrowvert f- \mathcal{I}_m f \arrowvert_{m+1}^2 . \label{Pyth}
\end{equation}

The approximation properties of Hermite interpolants of smooth functions are generally what one expects. Restricting to $d=1$ and using $\| \cdot \|$ to denote the $L^2$-norm we have, for mesh and function-independent constants, $C$ 
\begin{eqnarray}
\| f - \mathcal{I}_m f \| & \leq & C h^{2m+2} \| \frac {\partial^{2m+2} f}{\partial x^{2m+2}} \| , \label{L2interperr} \\
\| \frac {\partial^{m+1} (f- \mathcal{I}_m f )}{\partial x^{m+1} }\| & \leq & C h^{m+1} \| \frac {\partial^{2m+2} f}{\partial x^{2m+2}} \| . \label{Hm+1interperr}
\end{eqnarray}
We will also use (see \cite{goodrich2006hermite}):
\begin{equation}
\| f - \mathcal{I}_m \| \leq  C h^{m+1} \| \frac {\partial^{m+1} (f- \mathcal{I}_m f)}{\partial x^{m+1}} \| . \label{L2interperr2}
\end{equation}

\subsection{Conservation}
We now construct variables which are conserved during the solution process. We begin by introducing the operators $\mathcal{S}_{\pm}$ defined via Fourier series
\begin{equation}
    \widehat{\mathcal{S}_{\pm} f}=e^{\pm i c (k_1^2 + \ldots + k_d^2)^{1/2} \Delta t/2} \hat{f}(k) , \label{Sdef}
\end{equation}
where $c$ denotes the speed of wave propagation. Note that for $d=1$ the branch can be chosen so that these are simply shift operators
\begin{displaymath}
\mathcal{S}_{\pm} f(x) = f(x \pm c \Delta t/2), \ \  d=1,
\end{displaymath}
where $c$ is the speed of wave propagation. 
By Parseval's Theorem these operators preserve all Sobolev seminorms and in addition satisfy
\begin{equation}
    \mathcal{S}_{+} \mathcal{S}_{-}=\mathcal{S}_{-} \mathcal{S}_{+} = I . \label{Spminv}
\end{equation}
Using them we can write down a two-level conservation condition which holds at the continuous level.

\subsubsection{One space dimension}

For $d=1$, note that (\ref{eq:final1}) translated to Fourier variables leads to the exact solution formula
\begin{eqnarray*}
    \hat{p}(k,t+\Delta t) & = & \hat{p}(k,t) + 2 \sum_{j=1,{\rm odd}} \frac {(ik c \Delta t/2)^j}{j!} \hat{v}(k,t+\Delta t/2) \\ & = & \hat{p}(k,t) + \left( e^{ick \Delta t/2} - e^{-ick \Delta t/2} \right) \hat{v}(k,t+\Delta t/2) ,
\end{eqnarray*}
which, along with the analogous formulas for the evolution of $v$ implies
\begin{eqnarray}
    p(x,t+\Delta t) & = & p(x,t) + \mathcal{S}_{+} v(x,t+\Delta t/2) - \mathcal{S}_{-} v(x,t+\Delta t/2) , \label{p1evolve} \\
    v(x,t+\Delta t/2) & = & v(x,t-\Delta t/2) + \mathcal{S}_{+} p(x,t) - \mathcal{S}_{-} p(x,t) . \label{v1evolve} 
\end{eqnarray}
Define
\begin{eqnarray}
    P_{\pm} (x,t) & = & p(x,t) \mp \mathcal{S}_{\pm} v(x,t-\Delta t/2) , \label{Ppmdef} \\
    V_{\pm} (x,t+\Delta t/2) & = & v(x,t+\Delta t/2) \mp \mathcal{S}_{\pm} p(x,t) . \label{Vpmdef}
\end{eqnarray}
Then, using (\ref{Spminv}) we can rewrite (\ref{p1evolve})-(\ref{v1evolve}) as
\begin{eqnarray}
P_{\pm} (x,t+\Delta t) & = & \mp \mathcal{S}_{\mp} V_{\pm} (x,t+\Delta t/2) , \label{Ppmev1}\\
V_{\pm} (x,t+\Delta t /2) & = & \mp \mathcal{S}_{\mp} P_{\pm} (x, t) , \label{Vpmev1}
\end{eqnarray}
which yields the conservation conditions, valid for any Sobolev seminorm,
\begin{eqnarray}
    \| P_{\pm} (\cdot,t+\Delta t) \| & = & \| V_{\pm} (\cdot,t+\Delta t/2) \|, \label{cons1d_contP} \\
    \| V_{\pm} (\cdot,t+\Delta t/2) \| & = & \| P_{\pm} (\cdot,t) \| . \label{cons1d_contV}
\end{eqnarray}

Now consider the approximate evolution. The essential point is that the Hermite data updated according to (\ref{eq:final1}) is in fact the exact data for the update of the global piecewise polynomial Hermite interpolant so long as the CFL condition $c \Delta t/\Delta x < 1$ is satisfied; this simply follows from domain-of-dependence considerations. Therefore we can write the discrete evolution as
\begin{eqnarray}
    p^h(x,t+\Delta t) & = & p^h(x,t) + \mathcal{I}_m \mathcal{S}_{+} v^h(x,t+\Delta t/2) \nonumber 
    \\ & & - \mathcal{I}_m \mathcal{S}_{-} v^h(x,t+\Delta t/2) , \label{ph1evolve} \\
    v^h(x,t+\Delta t/2) & = & v^h(x,t-\Delta t/2) + \mathcal{I}_m \mathcal{S}_{+} p^h(x,t) \nonumber \\ & & - \mathcal{I}_m \mathcal{S}_{-} p^h(x,t) , \label{vh1evolve} 
\end{eqnarray}
Defining $P_{\pm}^h$, $V_{\pm}^h$ as in (\ref{Ppmdef})-(\ref{Vpmdef}) with the discrete variables $p^h$ and $v^h$ replacing the continuous ones and using the fact that $p^h=\mathcal{I}_m p^h$, $v^h = \mathcal{I}_m v^h$, we derive evolution equations
\begin{eqnarray}
    P_{\pm}^h (x,t+\Delta t) & = & \mp \mathcal{I}_m \mathcal{S}_{\mp} V_{\pm}^h (x,t+\Delta t/2) \nonumber \\ & & \mp (1-\mathcal{I}_m) \mathcal{S}_{\pm} V_{\mp}^h (x,t+\Delta t/2) , \label{Ppmhev1} \\
    V_{\pm}^h (x,t+\Delta t/2) & = & \mp \mathcal{I}_m \mathcal{S}_{\mp} P_{\pm}^h (x,t) \nonumber \\ & & \mp (1-\mathcal{I}_m) \mathcal{S}_{\pm} P_{\mp}^h (x,t) . \label{Vpmhev1} 
\end{eqnarray}
Now we get conservation conditions for the sums of the Sobolev seminorms of $P_{\pm}^h$ and $V_{\pm}^h$ using the orthogonality lemma (\ref{Pyth}). Define
\begin{eqnarray}
Q^h(t) & = & \left\arrowvert P_{+}^h (\cdot , t ) \right\arrowvert_{m+1}^2 + \left\arrowvert P_{-}^h (\cdot , t) \right\arrowvert_{m+1}^2 \label{Qdef} \\
R^h (t+\Delta t/2) & = & \left\arrowvert V_{+}^h (\cdot , t+\Delta t/2 ) \right\arrowvert_{m+1}^2 + \left\arrowvert V_{-}^h (\cdot , t+\Delta t/2) \right\arrowvert_{m+1}^2 . \label{Rdef}
\end{eqnarray}
Then 
\begin{eqnarray}
Q^h(t+\Delta t) & = & \left\arrowvert \mathcal{I}_m \mathcal{S}_{-} V_{+}^h (\cdot , t+ \Delta t/2) \right\arrowvert_{m+1}^2 + \left\arrowvert (1-\mathcal{I}_m) \mathcal{S}_{+} V_{-}^h (\cdot, t + \Delta t/2) \right\arrowvert_{m+1}^2 \nonumber \\ & & + \left\arrowvert \mathcal{I}_m \mathcal{S}_{+} V_{-}^h (\cdot , t+ \Delta t/2) \right\arrowvert_{m+1}^2 + \left\arrowvert (1-\mathcal{I}_m) \mathcal{S}_{-} V_{+}^h (\cdot, t + \Delta t/2) \right\arrowvert_{m+1}^2 \nonumber \\ & = & \left\arrowvert \mathcal{S}_{-} V_{+}^h (\cdot, t+\Delta t/2) \right\arrowvert_{m+1}^2 + \left\arrowvert \mathcal{S}_{+} V_{-}^h (\cdot, t+\Delta t/2) \right\arrowvert_{m+1}^2 \nonumber \\ & = & R^h (t+\Delta t/2) . \label{cons1d_hP}
\end{eqnarray}
Similarly
\begin{equation}
R^h (t+ \Delta t/2) = Q^h (t) , \label{cons1dh_V}
\end{equation} 
so we have for all $n$
\begin{equation}
Q^h( n \Delta t) = R^h((n+1/2) \Delta t) = Q^h(0) . \label{fullcons1d}
\end{equation}

\subsubsection{Extension to $d$ space dimensions}

In higher dimensions we consider the generalization of (\ref{eq:waveModel})
\begin{eqnarray}
\frac {\partial p}{\partial t} & = & c \nabla \cdot \mathbf{v} \label{wavend} \\
\frac {\partial \mathbf{v}}{\partial t} & = & c \nabla p , \nonumber
\end{eqnarray}
though other systems, such Maxwell's equations, could be similarly treated. The formulas analogous to (\ref{eq:final1}) are now

\begin{align}
p(\mathbf{x},t+\Delta t)  = 
p(\mathbf{x},t) + 2 \sum_{j = 1, odd}  \frac {( c \Delta t/2)^j}{j!}  \nabla^{j-1} \nabla \cdot \mathbf{v} ( \mathbf{x} , t+ \Delta t/2) 
, \label{pdup} \\
\mathbf{v} ( \mathbf{x} , t+ \Delta t/2)  =
\mathbf{v} (\mathbf{x}, t-\Delta t/2 ) +  
2 \sum_{j = 1, odd}  \frac{(c \Delta t/2)^j}{j!} \nabla^{j-1} \nabla  p( \mathbf{x} , t) . \label{vdup}
\end{align}

In Fourier space, denote the normalized $\mathbf{k}$-vector by
\begin{displaymath}
\hat{\mathbf{k}} = \frac {\mathbf{k}}{\arrowvert k \arrowvert},
\end{displaymath}
noting that in the formulas below the undefined value at $\mathbf{k}=\mathbf{0}$ plays no role. However we will need to specify its value when defining the conserved quantities and we do so by setting it to be some arbitrarily-chosen unit $d$-vector. We then find
\begin{eqnarray}
\hat{p}(\mathbf{k},t+\Delta t) & = & \hat{p}(\mathbf{k},t) \label{phatup} \\ & & + \left( e^{i c(k_1^2 + \ldots + k_d^2 )^{1/2} \Delta t/2 } - e^{-ic( k_1^2 + \ldots + k_d^2)^{1/2} \Delta t/2} \right) \hat{\mathbf{k}} \cdot \hat{\mathbf{v}} ( \mathbf{k} , t+ \Delta t/2) , \nonumber \\
\hat{\mathbf{v}} ( \mathbf{k} , t+ \Delta t/2) & = & \hat{\mathbf{v}} (\mathbf{k}, t-\Delta t/2 ) \label{vhatup} \\ & & +  \left( e^{ic (k_1^2 + \ldots + k_d^2)^{1/2} \Delta t/2 } - e^{-ic(k_1^2 + \ldots + k_d^2 )^{1/2}  \Delta t/2} \right) \hat{\mathbf{k}}  \hat{p} ( \mathbf{k} , t) . \nonumber  
\end{eqnarray}
Defining $\mathcal{G}$ by the operator given in Fourier space by multiplication by $\hat{\mathbf{k}}$ and
by $\mathcal{G}^{\ast}$ its adjoint, noting also that $\mathcal{G}^{\ast} \mathcal{G}=I$ we rewrite the continuous evolution formulas. We also define
\begin{displaymath}
\mathbf{v}_d = \mathcal{G} \mathcal{G}^{\ast} \mathbf{v}, \ \ \mathbf{v}_c = \mathbf{v} - \mathbf{v}_d 
\end{displaymath}
and note that
\begin{displaymath}
\mathcal{G}^{\ast} \mathbf{v}_d = \mathcal{G}^{\ast} \mathbf{v}, \ \ \mathcal{G}^{\ast} \mathbf{v}_c = 0.
\end{displaymath}
Then we have 
\begin{eqnarray}
p(\mathbf{x},t+\Delta t) & = & p(\mathbf{x},t) + \mathcal{S}_{+} \mathcal{G}^{\ast} \mathbf{v}_d (\mathbf{x},t+ \Delta t/2) \nonumber \\ & & - \mathcal{S}_{-} \mathcal{G}^{\ast} \mathbf{v}_d (\mathbf{x},t+\Delta t/2) , \label{pdopup} \\
\mathbf{v}_d (\mathbf{x},t+\Delta t/2) & = & \mathbf{v}_d (\mathbf{x},t-\Delta t/2) + \mathcal{G} \mathcal{S}_{+} p(\mathbf{x},t) \nonumber \\ & & - \mathcal{G} \mathcal{S}_{-} p(\mathbf{x},t) , \label{vddopup} \\
\mathbf{v}_c (\mathbf{x},t+\Delta t/2) & = & \mathbf{v}_c (\mathbf{x},t- \Delta t/2) . \label{vcdopup} 
\end{eqnarray}
Introducing 
\begin{eqnarray}
P_{\pm} (\mathbf{x},t) & = & p(\mathbf{x},t) \mp \mathcal{S}_{\pm} \mathcal{G}^{\ast} \mathbf{v}_d  (\mathbf{x},t-\Delta t/2) , \label{Pddef} \\
\mathbf{V}_{\pm} (\mathbf{x}, t+\Delta t/2) & = & \mathbf{v}_d (\mathbf{x},t) \mp \mathcal{G} \mathcal{S}_{\pm} p(\mathbf{x},t) , \label{Vddef}
\end{eqnarray}
we find
\begin{eqnarray}
P_{\pm} (\mathbf{x},t+\Delta t) & = & \mp \mathcal{S}_{\mp} \mathcal{G}^{\ast} \mathbf{V}_{\pm} (\mathbf{x},t+\Delta t/2) , \label{consdPup} \\ 
\mathbf{V}_{\pm} (\mathbf{x}, t+\Delta t/2 ) & = & \mp \mathcal{G} \mathcal{S}_{\mp} P_{\pm} (\mathbf{x},t) . \label{consdVup} 
\end{eqnarray}
Using (\ref{vcdopup}), (\ref{consdPup})-(\ref{consdVup}) and the fact that
\begin{displaymath}
\| \mathcal{G}^{\ast} \mathbf{V}_{\pm} \| = \| \mathbf{V}_{\pm} \| 
\end{displaymath}
the conservation of arbitrary Sobolev seminorms of $\mathbf{v}_c$, $P_{\pm}$, and $\mathbf{V}_{\pm}$ follows:
\begin{eqnarray}
\| P_{\pm} (\cdot ,t+\Delta t) \| & = & \| \mathbf{V}_{\pm} (\cdot, t+ \Delta t/2) \| \\
\| \mathbf{V}_{\pm} (\cdot ,t+\Delta t/2) \| & = & \| P_{\pm} (\cdot, t) \|  \\
\| \mathbf{v}_c (\cdot , t+\Delta t/2) \| & = & \| \mathbf{v}_c (\cdot , t- \Delta t/2) \| . 
\end{eqnarray}

Similarly we can establish conservation conditions for the Hermite-leapfrog updates. The main new issue which arises is the
fact that the interpolation operators do not commute with $\mathcal{G}$. Defining
\begin{displaymath}
\mathbf{v}_d^h = \mathcal{G} \mathcal{G}^{\ast} \mathbf{v}^h, \ \ \mathbf{v}_c^h = \mathbf{v}^h - \mathbf{v}_d^h,
\end{displaymath}
\begin{displaymath}
P_{\pm}^h = p^h \mp \mathcal{S}_{\pm} \mathcal{G}^{\ast} \mathbf{v}_d^h , \ \ \mathbf{V}_{\pm}^h = \mathbf{v}_d^h \mp \mathcal{G} \mathcal{S}_{\pm} p^h ,
\end{displaymath}
we derive the following update formulas
\begin{eqnarray}
P_{\pm}^h (\mathbf{x},t+\Delta t) & = & \mp \mathcal{I}_m \mathcal{S}_{\mp} \mathcal{G}^{\ast} \mathbf{V}_{\pm}^h (\mathbf{x},t+\Delta t/2) \label{Pdpmhev} \\ & & \mp (1-\mathcal{I}_m) \mathcal{S}_{\pm} \mathcal{G}^{\ast} \mathbf{V}_{\mp}^h (\mathbf{x},t+\Delta t/2) , \nonumber \\
\mathbf{V}_{\pm}^h (\mathbf{x}, t+\Delta t/2) + \mathbf{v}_c^h (\mathbf{x},t+\Delta t/2) & = & \mp \mathcal{I}_m \left( \mathcal{G} \mathcal{S}_{\mp} P_{\pm}^h (\mathbf{x},t) - \mathbf{v}_c^h(\mathbf{x},t-\Delta t/2) \right) \label{Vdmphev} \\ & & \mp (1-\mathcal{I}_m) \left( \mathcal{G} \mathcal{S}_{\pm} P_{\mp}^h (\mathbf{x},t) - \mathbf{v}_c^h (\mathbf{x},t-\Delta t/2) \right). \nonumber 
\end{eqnarray}
Define
\begin{eqnarray}
Q^h(t) & = & \left\arrowvert P_{+}^h(\cdot,t) \right\arrowvert_{m+1}^2 + \left\arrowvert P_{-}^h (\cdot, t) \right\arrowvert_{m+1}^2 + \left\arrowvert \mathbf{v}_c^h (\cdot,t-\Delta t/2) \right\arrowvert_{m+1}^2 \\ R^h (t+\Delta t/2) & = & \left\arrowvert \mathbf{V}_{+}^h (\cdot, t+\Delta t/2) \right\arrowvert_{m+1}^2 + \left\arrowvert \mathbf{V}_{-}^h (\cdot,t + \Delta t/2 ) \right\arrowvert_{m+1}^2 \\ & &  + \left\arrowvert \mathbf{v}_c^h (\cdot,t+\Delta t/2) \right\arrowvert_{m+1}^2 . \nonumber 
\end{eqnarray}
Then, since $\mathbf{v}_c^h$ is orthogonal in the Sobolev semi-inner products to the range of $\mathcal{G}$ we find as above
\begin{eqnarray}
Q^h(t+\Delta t) & = & R^h (t+\Delta t/2) , \\ R^h(t+\Delta t/2)  & = & Q^h(t) ,
\end{eqnarray}
from which we finally conclude 
\begin{equation}
Q^h (n \Delta t) = R^h ((n+1/2) \Delta t) = Q^h (0) .  \label{fullconsd}
\end{equation}

\subsection{Convergence in one dimension}

We now exploit the energy estimates along with the approximation properties of Hermite interpolation (\ref{L2interperr})-(\ref{L2interperr2}) to establish convergence. We define the errors by
\begin{displaymath}
e_p = p-p^h, \ \ e_v = v-v^h .
\end{displaymath}
They satisfy the evolution formula
\begin{eqnarray}
e_p(x,t+\Delta t) & = & e_p (x,t) + \mathcal{I}_m \mathcal{S}_{+} e_v (x,t+\Delta t/2) \label{epev} \\
& & - \mathcal{I}_m \mathcal{S}_{-} e_v (x, t+\Delta t/2) + (1-\mathcal{I}_m) p(x,t+\Delta t) , \nonumber \\
e_v(x,t+\Delta t/2) & = & e_v (x,t-\Delta t/2) + \mathcal{I}_m \mathcal{S}_{+} e_p (x,t) \label{evev} \\
& & - \mathcal{I}_m \mathcal{S}_{-} e_p (x, t) + (1-\mathcal{I}_m) v(x,t+\Delta t/2) . \nonumber 
\end{eqnarray}
Introducing the variables corresponding to the conserved quantities
\begin{eqnarray}
    E_{p,\pm} (x,t) & = & e_p(x,t) \mp \mathcal{S}_{\pm} e_v(x,t-\Delta t/2) , \label{EPpmdef} \\
    E_{v,\pm} (x,t+\Delta t/2) & = & e_v(x,t+\Delta t/2) \mp \mathcal{S}_{\pm} e_p(x,t) , \label{EVpmdef}
\end{eqnarray}
and following the previous calculations we derive 
\begin{eqnarray}
    E_{p,\pm} (x,t+\Delta t) & = & \mp \mathcal{I}_m \mathcal{S}_{\mp} E_{v,\pm} (x,t+\Delta t/2) \nonumber \\ & & \mp (1-\mathcal{I}_m) \mathcal{S}_{\pm} E_{v,\mp} (x,t+\Delta t/2) \nonumber \\ & & 
    + (1-\mathcal{I}_m) (p(x,t+\Delta t)-p(x,t)) , \label{EPpmev1} \\
    E_{v,\pm} (x,t+\Delta t/2) & = & \mp \mathcal{I}_m \mathcal{S}_{\mp} E_{p,\pm} (x,t) \nonumber \\ & & \mp (1-\mathcal{I}_m) \mathcal{S}_{\pm} E_{p,\mp} (x,t) \nonumber \\ & & + (1-\mathcal{I}_m) (v(x,t+\Delta t/2) - v(x,t-\Delta t/2)) . \label{EVpmev1} 
\end{eqnarray}
Noting that for sufficiently smooth solutions 
\begin{eqnarray*}
\left\arrowvert (1-\mathcal{I}_m) (p(x,t+\Delta t)-p(x,t)) \right\arrowvert_{m+1} & = & O\left(\Delta t \cdot h^{m+1} \right), \\ \left\arrowvert (v(x,t+\Delta t/2) - v(x,t-\Delta t/2)) \right\arrowvert_{m+1} & = &  O\left(\Delta t \cdot h^{m+1} \right), 
\end{eqnarray*}
and defining
\begin{eqnarray*}
\mathcal{E}_p (t) & = & \left\arrowvert E_{p,+}(\cdot,t) \right\arrowvert_{m+1}^2 + \left\arrowvert E_{p,-} (\cdot,t) \right\arrowvert_{m+1}^2 , \\
\mathcal{E}_v (t+\Delta t/2) & = & \left\arrowvert E_{v,+}(\cdot,t+\Delta t/2) \right\arrowvert_{m+1}^2 + \left\arrowvert E_{v,-} (\cdot,t+\Delta t/2) \right\arrowvert_{m+1}^2 , 
\end{eqnarray*}
we find, where here and in what follows $C$ denotes a mesh-independent quantity which will depend on the exact solution,
\begin{eqnarray}
\mathcal{E}_p(t + \Delta t) & \leq & \mathcal{E}_v (t+\Delta t/2) + C \Delta t \cdot h^{m+1} \sqrt{\mathcal{E}_v(t+\Delta t/2)} \nonumber \\ & & + C \Delta t^2 h^{2m+2} , \label{Epineq} \\ 
\mathcal{E}_v(t + \Delta t/2) & \leq & \mathcal{E}_p (t) + C \Delta t \cdot h^{m+1} \sqrt{\mathcal{E}_p(t)} + C \Delta t^2 h^{2m+2} . \label{Evineq}
\end{eqnarray}
Defining
\begin{displaymath}
\hat{\mathcal{E}}^n = \max_{j \leq n} \left(\mathcal{E}_p(j \Delta t), \mathcal{E}_n( (j-1/2) \Delta t) \right) 
\end{displaymath}
we have by summing (\ref{Epineq})-(\ref{Evineq}) over $j$
\begin{displaymath}
\hat{\mathcal{E}}^n \leq \hat{\mathcal{E}}^0 + Ct_n h^{m+1} \sqrt{\hat{\mathcal{E}}^n} + Ct_n \Delta t \cdot h^{2m+2} . 
\end{displaymath}
Assuming, as would be true for smooth initial data and initialization of the scheme by Hermite interpolation 
\begin{displaymath}
\hat{\mathcal{E}}^0 \leq C h^{2m+2},
\end{displaymath}
we derive our first error estimate
\begin{equation}
\hat{\mathcal{E}}^n \leq C(1+t_n^2) h^{2m+2} . \label{errest1}
\end{equation}

Now return to (\ref{EPpmev1})-(\ref{EVpmev1}), which we rewrite as
\begin{eqnarray*}
E_{p,\pm} (x,t+\Delta t) & = & \mp \mathcal{S}_{\mp} E_{v,\pm} (x,t+\Delta t/2) \nonumber \\ & & \mp (1-\mathcal{I}_m) \left( \mathcal{S}_{\pm} E_{v,\mp} (x,t+\Delta t/2) - \mathcal{S}_{\mp} E_{v,\pm} (x,t+\Delta t/2) \right)  \\ & & 
    + (1-\mathcal{I}_m) (p(x,t+\Delta t)-p(x,t)) ,  \\
    E_{v,\pm} (x,t+\Delta t/2) & = & \mp \mathcal{S}_{\mp} E_{p,\pm} (x,t) \nonumber \\ & & \mp (1-\mathcal{I}_m) \left( \mathcal{S}_{\pm} E_{p,\mp} (x,t) - \mathcal{S}_{\mp} E_{p,\pm} (x,t) \right) \\ & & + (1-\mathcal{I}_m) (v(x,t+\Delta t/2) - v(x,t-\Delta t/2)) . 
\end{eqnarray*}
Computing the $L^2$ norms and using (\ref{L2interperr2}) along with (\ref{errest1}) we deduce
\begin{eqnarray*}
\| E_{p,\pm} (\cdot,t+ \Delta t) \| & \leq & \| E_{v,\pm} (\cdot,t+\Delta t/2) \| + C (1+t) h^{2m+2} , \\
\| E_{v,\pm} (\cdot,t+ \Delta t/2) \| & \leq & \| E_{p,\pm} (\cdot,t) \| + C (1+t) h^{2m+2} .
\end{eqnarray*}
Assuming sufficiently accurate initial data and in addition that for some $\eta > 0$, $\Delta t > \eta h$ we
sum these inequalities and deduce our second error estimate:
\begin{eqnarray}
\| E_{p,\pm} (\cdot,n \Delta t) \| & \leq &   C (1+t_n^2) h^{2m+1} , \label{Eperrest} \\
\| E_{v,\pm} (\cdot,(n-1/2) \Delta t) \| & \leq &  C (1+t_n^2) h^{2m+1} . \label{Everrest}
\end{eqnarray}

Finally, we rewrite (\ref{EPpmdef})-(\ref{EVpmdef})  
\begin{eqnarray*}
 e_p(x,t) & = &   E_{p,\pm} \pm \mathcal{S}_{\pm} e_v(x,t-\Delta t/2) ,  \\
e_v(x,t+\Delta t/2) & = & E_{v,\pm} (x,t+\Delta t/2) \pm \mathcal{S}_{\pm} e_p(x,t) , 
\end{eqnarray*}
which yields after taking norms
\begin{eqnarray*}
 \| e_p(\cdot,t) \| & \leq &   \| e_v(\cdot,t-\Delta t/2) \| + C(1+t^2) h^{2m+1},  \\
\| e_v(\cdot, t+\Delta t/2) & = & \| e_p(\cdot,t) \| + C(1+t^2) h^{2m+1} . 
\end{eqnarray*}
Summing this we get an estimate of the solution error:
\begin{eqnarray}
\| e_p (\cdot,n \Delta t) \| & \leq &   C (1+t_n^3) h^{2m} , \label{eperrest} \\
\| e_v (\cdot,(n-1/2) \Delta t) \| & \leq &  C (1+t_n^3) h^{2m} . \label{everrest}
\end{eqnarray}

We note that in the numerical experiments we do not see the cubic growth in $t$ which appears in (\ref{eperrest})-(\ref{everrest}). There is a technical barrier to the extension of the error analysis to higher spatial dimensions, namely that the seminorm used in the stability analysis cannot control the full interpolation error due to the fact that it is zero for functions of fewer than $d$ variables. However we have not noticed any degradation in convergence rates for our experiments with $d=2,3$. What the analysis misses is the fact that, for $m$ even, we observe convergence at order $2m+2$ rather than $2m$. In the next section we study the dispersion relation for the method which suggests an explanation for this phenomenon. 

\subsection{Dispersion}

As a complement to the stability and error analysis given above, we also consider the dispersive properties of the method. Here we again assume spatial periodicity and expand the discrete solution data in a Fourier series:
\begin{eqnarray*}
\tilde{P}^h (x_j,t) & = & \sum_k \hat{\tilde{P}} (k,t) e^{i k x_j} , \\
\tilde{V}^h (x_{j+1/2},t) & = & \sum_k \hat{\tilde{V}} (k,t) e^{i k x_{j+1/2}} .
\end{eqnarray*}
Here $\tilde{P}^h$ and $\tilde{V}^h$ denote the $m+1$-vectors of approximate function and derivative data at the nodes so that $\hat{\tilde{P}}$ and $\hat{\tilde{V}}$ are also $m+1$-vectors. Written in these variables the evolution takes the form:
\begin{eqnarray}
\hat{\tilde{P}}(k,t+\Delta t) & = & \hat{\tilde{P}}(k,t) + \hat{D}(k,\Delta t) \hat{\tilde{V}} (k,t+\Delta t/2) , \label{Ptildehatev} \\
\hat{\tilde{V}}(k,t+\Delta t/2) & = & \hat{\tilde{V}}(k,t-\Delta t/2) + \hat{D}(k,\Delta t) \hat{\tilde{P}} (k,t) , \label{Vtildehatev} 
\end{eqnarray}
where the existence of matrix $\hat{D}$ follows from the translation-invariance of the method. Obviously the solutions can be expressed using the eigenvalues, $\kappa$, and eigenvectors, $\hat{w}$ of $\hat{D}$ which we will compare to
the exact values, 
\begin{displaymath}
\kappa_c=2i \sin{(ck\Delta t/2)}, \ \ \hat{w}_c=\left( 1 \ ik \ \ldots \ (ik)^m/m! \right)^T .  
\end{displaymath}
To make this comparison at leading order we rescale the problem so that $h=1$ and
$\Delta t=\lambda$ and make an expansion for $k \ll 1$. Using the accuracy properties of the Hermite interpolation operators we have
\begin{equation}
\hat{D} \hat{w}_c = \kappa_c \hat{w}_c + O(k^{2m+2}). \label{Dtrunc}
\end{equation}
We thus deduce the perturbation equation
\begin{equation}
\hat{D}(\hat{w}-\hat{w}_c) = (\kappa - \kappa_c) \hat{w}_c + \kappa (\hat{w}-\hat{w}_c) + O(k^{2m+2}). 
\label{Dpert}
\end{equation}
To leading order we may replace $\hat{D}(k)$ by $\hat{D}(0)$, $\hat{w}_c(k)$ by $\hat{w}_c(0)=e_1$, and, since
$\kappa_c(0)=0$, ignore $\kappa (\hat{w}-\hat{w}_c)$. We then have two possibilities:
\begin{description}
\item[i.] $0$ is a simple eigenvalue of $\hat{D}(0)$, in which case there is a solution to (\ref{Dpert}) satisfying
\begin{equation}
\kappa-\kappa_c = O(k^{2m+2}), \ \ \hat{w}-\hat{w}_c = O(k^{2m+2}) , \label{goodpert} 
\end{equation}
\item[ii.] 0 is not a simple eigenvalue of $\hat{D}(0)$ in which case we do expect (\ref{goodpert}) to hold. 
\end{description}
To apply this result we check for solutions for a second eigenvector corresponding to the eigenvalue $0$. 
\begin{equation}
\hat{D}(0) \hat{g} = 0 , \label{ev2}
\end{equation}
or a generalized eigenvector
\begin{equation}
\hat{D}(0) \hat{g} = e_1 . \label{gev}
\end{equation}
Note that solutions of these equations correspond to the existence of a nonconstant polynomial $g(x)$ of degree $2m+1$ satisfying for $j=0, \ldots ,m$
\begin{displaymath}
\frac {d^j g}{dx^j} (1/2) - \frac {d^j g}{dx^j} (-1/2) = \frac {d^j g}{dx^j} (\lambda/2) - \frac {d^j g}{dx^j} (-\lambda/2) =0 ,
\end{displaymath}
or the above equations except
\begin{displaymath}
g(\lambda/2)-g(-\lambda/2)=1.
\end{displaymath}
Moreover, since a constant function solves the homogeneous equation, we can in addition assume that $g(0)=0$. Expanding $g$
\begin{displaymath}
g=\sum_{j=1}^{2m+1} g_j x^j 
\end{displaymath}
we see that they represent $2m+2$ linear equations with $2m+1$ unknown coefficients. Moreover, due to symmetry the equations involving even order derivatives only involve the $m+1$ coefficients of the odd powers of $x$ while equations involving odd order derivatives only involve the $m$ coefficients of even powers. The size of these subsystems differs in the cases $m=2n$ and $m=2n+1$:
\begin{description}
\item[$m=2n$:] We have $2n+2$ equations for $2n+1$ coefficients of odd powers and $2n$ equations for coefficients of even powers. As we have checked numerically, these are only solvable by $0$ when all equations are homogeneous. Thus $0$ is a simple eigenvalue and we have a solution which is accurate to
order $2m+2$.
\item[$m=2n+1$:] We have $2n+2$ equations for $2n+2$ coefficients of odd powers and $2n+2$ homogeneous equations for $2n+1$ coefficients of even powers. Although the latter equations are overdetermined, they are always solved by $0$. Numerically we have found that the equations for the coefficients of the odd powers is invertible; thus there is no second eigenvalue solving (\ref{ev2}) but there is a generalized eigenvector solving (\ref{gev}). Thus we do not in general expect a solution which is accurate to order $2m+2$, but only $2m$ as proven above.
\end{description}

\section{Numerical experiments in one-dimension}
\label{sec:NEx1D} 
To verify the observations in the previous section, we present numerical experiments which assess the accuracy and performance of the Hermite-leapfrog method. To compliment the experiments, we present comparisons with the classic Hermite method
(Dual Hermite method as referred to in \cite{vargas2017variations}). As a model equation we chose the following wave system
\begin{align} \label{eq:waveEq1D}
\frac{\partial p}{\partial t} &=  -c^2(x) \frac{\partial v}{\partial x} + z(x,t),  \quad\frac{\partial v}{\partial t} = -\frac{\partial p}{\partial x} \\ \nonumber
 p(x,t_0) &= f(x), \quad v(x,t_{0}+{\Delta t}/{2}) = g(x).\nonumber
\end{align}
Here $c(x)$ denotes the speed of the propagating wave and $z(x,t)$ is introduced to enforce a desired solution.
The time step, $\Delta t$ is chosen by introducing a CFL constant $C_{CFL} \in (0,1)$, such that
\[
\Delta t = 2 C_{CFL}\frac{h}{c_{max}}.
\] 
Choosing $C_{CFL}$ closer to one denotes a larger time step while choosing $C_{CFL}$ closer to zero corresponds to taking a smaller time-step. The variable $h$ corresponds to the length of the interpolation interval, and $c_{max}$ denotes the maximum speed of the propagating wave.

\subsection{Standing wave solution}
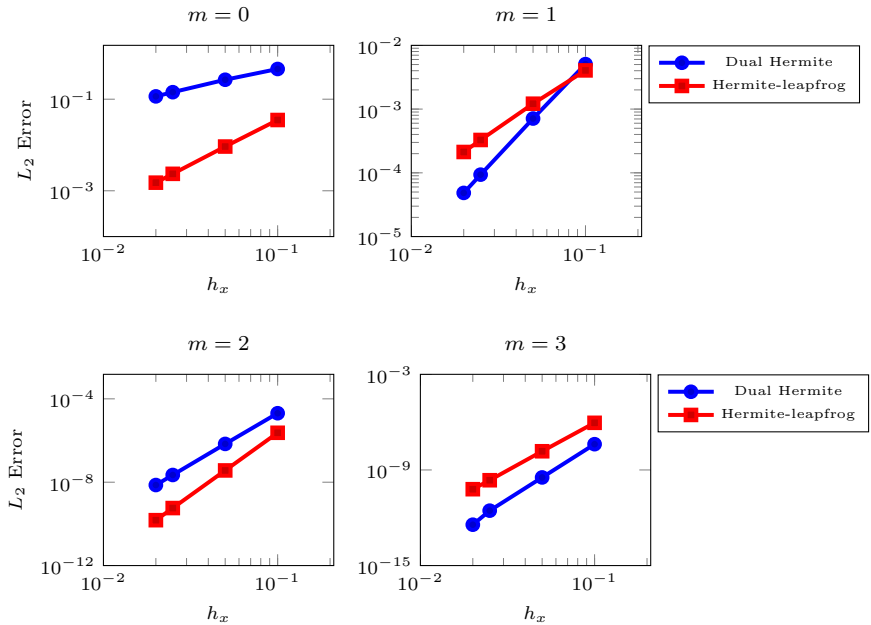
\begin{figure}[h!]
    \centering
\begin{subfloat}{
\begin{tikzpicture}
\begin{loglogaxis}[	
    legend style={font=\tiny},
	scale=0.36, 
    title={$m=0$},
    xlabel={$h_x$},
    ylabel={$L_{2}$ Error},
    xmin=0.01,xmax=0.21,
    ymin=0.0001,ymax=1.5,
    every axis plot/.append style={line width= 1.5pt}
]
\addplot table {hPlusRatesV2_T4_13/dualRateN0_C9.dat};
\addplot table {hPlusRatesV2_T4_13/hPlusRateN0_C9.dat};
\end{loglogaxis}
\end{tikzpicture}}
\end{subfloat}
\begin{subfloat}{
\begin{tikzpicture}
\begin{loglogaxis}[	
    legend style={font=\tiny},
	scale=0.36, 
    title={$m=1$},
    xlabel={$h_x$},
    xmin=0.01,xmax=0.21,
    ymin=1e-5,ymax=1e-2,
    legend pos= outer north east,
    every axis plot/.append style={line width= 1.5pt},
    legend entries={Dual Hermite, Hermite-leapfrog},
]
\addplot table {hPlusRatesV2_T4_13/dualRateN1_C9.dat};
\addplot table {hPlusRatesV2_T4_13/hPlusRateN1_C9.dat};
\end{loglogaxis}
\end{tikzpicture}}
\end{subfloat}
\begin{subfloat}{
\begin{tikzpicture}
\begin{loglogaxis}[	
    legend style={font=\tiny},
	scale=0.36, 
    title={$m=2$},
    xlabel={$h_x$},
    ylabel={$L_{2}$ Error},
    xmin=0.01,xmax=0.21,
    ymin=1e-12,ymax=1.5e-3,
    every axis plot/.append style={line width= 1.5pt}
]
\addplot table {hPlusRatesV2_T4_13/dualRateN2_C9.dat};
\addplot table {hPlusRatesV2_T4_13/hPlusRateN2_C9.dat};
\end{loglogaxis}
\end{tikzpicture}}
\end{subfloat}
\begin{subfloat}{
\begin{tikzpicture}
\begin{loglogaxis}[	
    legend style={font=\tiny},
	scale=0.36, 
    title={$m=3$},
    xlabel={$h_x$},
    xmin=0.01,xmax=0.21,
    ymin=1e-15, ymax=1e-3,
    legend pos= outer north east,
    every axis plot/.append style={line width= 1.5pt},
    legend entries={Dual Hermite, Hermite-leapfrog},
]
\addplot table {hPlusRatesV2_T4_13/dualRateN3_C9.dat};
\addplot table {hPlusRatesV2_T4_13/hPlusRateN3_C9.dat};
\end{loglogaxis}
\end{tikzpicture}}
\end{subfloat}
\caption{$L^2$ errors of an $m$ order Hermite-leapfrog and Dual Hermite methods when applied to the one-dimensional pressure-velocity system. The Hermite-leapfrog scheme provides a better approximation to the classic Hermite method when $m$ is chosen to be even. For these experiments the CFL constant is set to $CFL=0.9$ across all orders.}
\label{fig:hermitePlus1D}. 
\end{figure}

As a first example we consider propagation of a standing wave with unit wave speed. The analytic solution is chosen to be 
\[
p(x,t) = \cos(2\pi t) \sin(2\pi x).
\]
For these experiments the computational domain is defined to be the bi-unit interval and the solution is propagated to a final time of $T=4.13$. Figure \ref{fig:hermitePlus1D} reports the accuracy in the $L^2$ norm while Table \ref{table:hermitePlus1D} reports the observed rates of convergence. Numerically, we observe that even order Hermite-leapfrog schemes converge at rates of $O(h^{2m+2})$, as suggested by the analysis of the dispersion relation, while variation is observed if the method is of odd order. Furthermore, it can be observed that a Hermite-leapfrog scheme of even order provides a better approximation than the Dual Hermite method. 

\begin{table}[h!]
\centering
\scalebox{0.8}{   
\begin{tabular}{|l|l|l|l|l|l|l|l|l|l|l|l|l|}
\hline
             & \multicolumn{4}{l|}{$C_{CFL} = 0.1$} & \multicolumn{4}{l|}{$C_{CFL} = 0.5$} & \multicolumn{4}{l|}{$C_{CFL} = 0.9$} \\ \hline
Order - $m$            & 0     & 1    & 2    & 3    & 0     & 1     & 2    & 3    & 0     & 1    & 2    & 3    \\ \hline
Hermite-leapfrog & 1.90 & 1.96 & 5.97 & 5.65 & 1.94  & 1.91  & 5.98 & 7.09 & 2.01 & 2.02 & 6.00 & 5.87 \\ \hline
Classic Hermite & -  & 2.58 & 4.90 & 6.84 & 0.10  & 2.78  & 4.83 & 6.95 & 0.46 & 2.78 & 4.83 & 6.88 \\ \hline
\end{tabular}
}
\caption{Observed $L^2$ rates of convergence for the Hermite-leapfrog scheme and classic Hermite methods of $m^{th}$ order with varying CFL constants.}
\label{table:hermitePlus1D}
\end{table}

\subsection{Smoothly varying wave speed}
\begin{figure}[h!]
\centering
\begin{subfloat}{
\begin{tikzpicture}
\begin{loglogaxis}[	
    legend style={font=\tiny},
	scale=0.36, 
    title={$m=0$},
    xlabel={$h_x$},
    ylabel={$L_{2}$ Error},
    xmin=0.001,xmax=1,
    ymin=0.001,ymax=1,
    every axis plot/.append style={line width= 1.5pt},
]
\addplot table {LFRates2D/Dual_N0_L2.dat};
\addplot table {LFRates2D/HY_N0_L2.dat};
\end{loglogaxis}
\end{tikzpicture}}
\end{subfloat}
\begin{subfloat}{
\begin{tikzpicture}
\begin{loglogaxis}[	
    legend style={font=\tiny},
	scale=0.36, 
    title={$m=1$},
    xlabel={$h_x$},
    xmin=0.001,xmax=1,
    ymin=0.00001,ymax=1,
    legend pos= outer north east,
    every axis plot/.append style={line width= 1.5pt},
    legend entries={Dual Hermite, Hermite-leapfrog},
]
\addplot table {LFRates2D/Dual_N1_L2.dat};
\addplot table {LFRates2D/HY_N1_L2.dat};
\end{loglogaxis}
\end{tikzpicture}}
\end{subfloat}
\begin{subfloat}{
\begin{tikzpicture}
\begin{loglogaxis}[	
    legend style={font=\tiny},
	scale=0.36, 
    title={$m=2$},
    xlabel={$h_x$},
    ylabel={$L_{2}$ Error},
    xmin=0.001,xmax=1,
    ymin=0.000000001,ymax=1,
    every axis plot/.append style={line width= 1.5pt},
]
\addplot table {LFRates2D/Dual_N2_L2.dat};
\addplot table {LFRates2D/HY_N2_L2.dat};
\end{loglogaxis}
\end{tikzpicture}}
\end{subfloat}
\begin{subfloat}{
\begin{tikzpicture}
\begin{loglogaxis}[	
    legend style={font=\tiny},
	scale=0.36, 
    title={$m=3$},
    xlabel={$h_x$},
    xmin=0.001,xmax=1,
    ymin=0.00000000001,ymax=1,
    legend pos= outer north east,
    every axis plot/.append style={line width= 1.5pt},
    legend entries={Dual Hermite, Hermite-leapfrog},
]
\addplot table {LFRates2D/Dual_N3_L2.dat};
\addplot table {LFRates2D/HY_N3_L2.dat};
\end{loglogaxis}
\end{tikzpicture}}
\end{subfloat}
\caption{The $L^2$ errors of various $m$ order Hermite-leapfrog and Dual Hermite methods when applied to the one-dimensional pressure-velocity system with spatially varying wave speed. The Hermite-leapfrog scheme provides a better approximation than the Dual-Hermite method when $m$ is chosen to be even. Here the CFL constant is set to $C_{CFL}=0.9$.}
\label{fig:varWavef}. 
\end{figure}
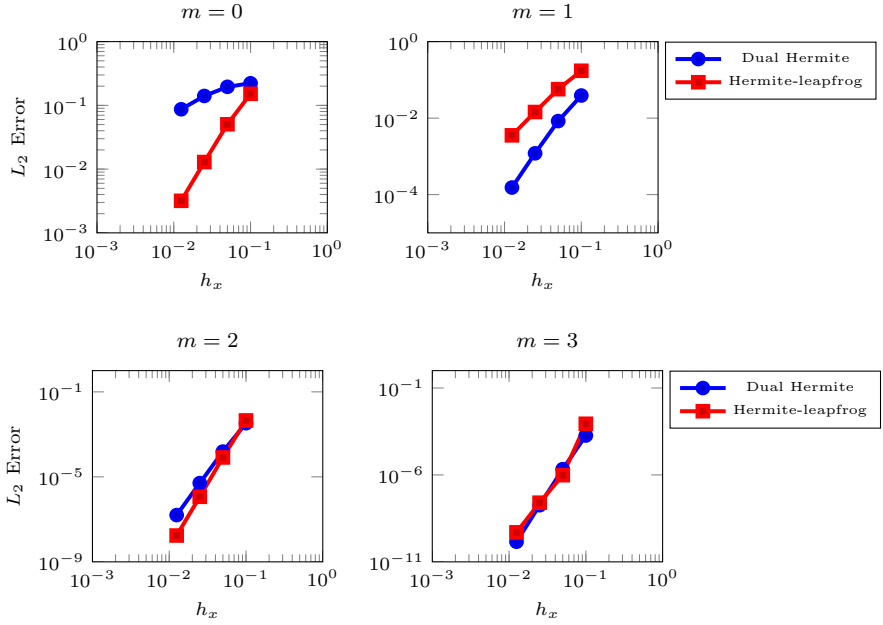

As a second set of numerical experiments we consider a spatially varying wave speed. To simplify establishing the Cauchy-Kovalesky recurrence relation we take a pre-processing step and expand the coefficient at each grid point via a local Taylor expansion
\[
c^2(x_i) = 1 + {\sin \left( x_i \right)}/{2} \approx \sum_{i=0}^{2m+1} \tilde{c}_i \left( \frac{x-x_i}{h} \right), \\ 
\]
where
\[
\tilde{c}_i = \frac{h!}{i!} \frac{\partial^i c^2(x)}{\partial x^i}. 
\]
For these experiments the domain is chosen to be $[0, 2\pi]$ and the solution is propagated to a final time of $T=3.2$. The analytic solution is chosen to be 
\begin{align} 
p(x,t) = \sin(x-t). 
\end{align}
Figure \ref{fig:varWavef} reports the observed accuracy in the $L^2$ norm while Table \ref{table:varWavef} reports the observed rates of convergence. These numerical experiments suggest rates of $O(h^{2m+2})$ for even $m$ and $O(h^{2m})$ for odd $m$.

\begin{table}[h!]
\centering
\scalebox{0.8}{   
\begin{tabular}{|l|l|l|l|l|l|l|l|l|l|l|l|l|}
\hline
             & \multicolumn{4}{l|}{$C_{CFL} = 0.1$} & \multicolumn{4}{l|}{$C_{CFL} = 0.5$} & \multicolumn{4}{l|}{$C_{CFL} = 0.9$} \\ \hline
$N$            & 0     & 1    & 2    & 3    & 0     & 1     & 2    & 3    & 0     & 1    & 2    & 3    \\ \hline
Hermite-leapfrog & 1.98 & 1.97 & 5.98 & 5.94 & 1.98 & 1.97 & 5.99 & 5.88 & 1.98 & 1.98 & 6.02 & 5.81 \\ \hline
Dual-Hermite     & 0.18 & 2.97 & 4.99 & 6.84 & 0.70 & 2.99 & 5.01 & 7.01 & 0.88 & 2.99 & 5.00 & 6.99 \\ \hline
\end{tabular}
}
\caption{Observed $L^2$ rates of convergence for various $m$ order Hermite-leapfrog and Dual Hermite methods when applied to the pressure-velocity system with smoothly varying coefficients. The numerical experiments are carried out for various $CFL$ constants.}
\label{table:varWavef}
\end{table}

\subsection{A note on improving rates of convergence}
The peculiar convergence behavior of the Hermite-leapfrog scheme has led us to investigate if modifications can be made which may yield more consistent convergence results.  In particular we have found a modification to the method which leads to consistent $O(h^{2m+2})$ convergence rates. The derivation of the modified Hermite-leapfrog scheme originates from applying the scheme to the following advection equation,
\begin{align}
u_t = u_x, \label{eq:adevec} 
\end{align}
and tracing out local truncation errors. The analog of subtracting Equations \ref{eq:option1} and \ref{eq:option2} for the advection equation yields a time-stepping scheme with local temporal truncation errors
\begin{align}
\frac{\partial^{r} u\left(x, t + \Delta t \right)}{\partial x^{r}} &= \frac{\partial^r u\left(x,t \right)}{\partial x^r} +  \sum_{i=1+r, odd}^{2m+1} c_{i-r} \frac{\partial^{i} u \left(x, t + {\Delta t}/{2} \right)}{\partial x^{i}} + \bigo(\Delta t^{2m+3-r}), \label{eq:hy1} \\
\frac{\partial^{l} u \left(x, t + \Delta t \right)}{\partial x^{l}} &= \frac{\partial^l u \left(x, t\right)}{\partial x^l} +  \sum_{i=1+l, even}^{2m} c_{i-l} \frac{\partial^{i} u\left(x, t + {\Delta t}/{2} \right)}{\partial x^{i}} + \bigo(\Delta t^{2m+2-l}) \label{eq:hy2}
\end{align}
where it is important to observe that Equation \ref{eq:hy1} is valid for even order derivatives, $r \in \{ 0,2,\dots$ \},  and Equation \ref{eq:hy2} is valid for odd order derivatives, $l \in \{ 1,3,\dots$ \}. Here $c_j$ are the coefficients from the Taylor series expansion. Drawing from the work of Appel\"{o} and co-authors in \cite{appelo2018hermite}, we construct a time-stepping scheme by subtracting expansions \ref{eq:option1} and \ref{eq:option2} leading to
\begin{align}
\frac{\partial^{r} u\left(x, t + \Delta t \right)}{\partial x^{r}} &= -\frac{\partial^r u (x, t) }{\partial x^r} +  \sum_{i=2+r}^{2m} c_{i-r} \frac{\partial^{i} u (x, t + {\Delta t}/{2})}{\partial x^{i}} + \bigo(\Delta t^{2m+2-r}), \label{eq:lfw1} \\
\frac{\partial^{l} u \left(x, t + \Delta t \right)}{\partial x^{l}} &= -\frac{\partial^l u \left(x, t \right)}{\partial x^l} +  \sum_{i=2+l}^{2m-1} c_{i-l} \frac{\partial^{l} u (x, t + {\Delta t}/{2})}{\partial x^{l}} + \bigo(\Delta t^{2m+3-l}), \label{eq:lfw2}
\end{align} 
where the indices remain to be $r=0,2,\dots$ and  $l=1,3,\dots$. As before, Equation \ref{eq:lfw1} is the local truncation error for even order derivatives, and \ref{eq:lfw2} is the local truncation for odd order derivatives. 

Spatial errors may be quantified by noting that Hermite-Birkhoff reconstruction leads to truncation errors of the form
\begin{align}
\frac{\partial u^r(x,t)}{\partial x^r} = 
\frac{\partial \tilde{u}^r(x,t)}{\partial x^r} + \bigo(h^{2m+2-r}) \quad r&=0,2,4,\dots \\
\frac{\partial u^l (x,t)}{\partial x^l} = 
\frac{\partial \tilde{u}^l(x,t)}{\partial x^l} + \bigo(h^{2m+3-l}) \quad l&=1,3,5,\dots. 
\end{align}
Here it is important to note that subsequent derivatives are approximated with the same order. 

By using the time-stepping scheme given by Equation \ref{eq:hy1} for even order derivatives and Equation \ref{eq:lfw2} for odd order derivatives we are choosing the best approximation within the framework. An immediate area of research would be to carryout a full truncation analysis, and a similar eigenvalue analysis to confirm $2m+2$ convergence rates. In the next section we present numerical evidence on stability and provide numerical estimates convergence rates for the modified Hermite-leapfrog method. 

\subsubsection{Modified Hermite-leapfrog: Advection Equation} 
\begin{figure}[h!]
\centering
\begin{subfloat}{
\begin{tikzpicture}
\begin{loglogaxis}[	
    legend style={font=\tiny},
    scale=0.37,
    title={$C_{CFL}=0.1$},
    xlabel={$h_x$},
    ylabel={$L_{2}$ Error},
    xmin=1.8e-2,xmax=0.4,
    ymin=1e-14,ymax=1.5,
    every axis plot/.append style={line width= 1.5pt},
]
\addplot table {advecPP/advecN1_CLF_1.dat};
\addplot table {advecPP/advecN2_CLF_1.dat};
\addplot table {advecPP/advecN3_CLF_1.dat};
\end{loglogaxis}
\end{tikzpicture}}
\end{subfloat}
\begin{subfloat}{
\begin{tikzpicture}
\begin{loglogaxis}[	
    legend style={font=\tiny},
    scale=0.37,
    title={$C_{CFL}=0.5$},
    xlabel={$h_x$},
   xmin=1.8e-2,xmax=0.4,,
    ymin=1e-14,ymax=1.5,
    every axis plot/.append style={line width= 1.5pt},
]
\addplot table {advecPP/advecN1_CLF_5.dat};
\addplot table {advecPP/advecN2_CLF_5.dat};
\addplot table {advecPP/advecN3_CLF_5.dat};
\end{loglogaxis}
\end{tikzpicture}}
\end{subfloat}
\begin{subfloat}{
\begin{tikzpicture}
\begin{loglogaxis}[	
    legend style={font=\tiny},
    scale=0.37,
    title={$C_{CFL}=0.9$},
    xlabel={$h_x$},
    xmin=1.8e-2,xmax=0.4,
    ymin=1e-14,ymax=1.5,
    legend pos= outer north east,
    every axis plot/.append style={line width= 1.5pt},
    legend entries={$m=1$,$m=2$,$m=3$},
]
\addplot table {advecPP/advecN1_CLF_9.dat};
\addplot table {advecPP/advecN2_CLF_9.dat};
\addplot table {advecPP/advecN3_CLF_9.dat};
\end{loglogaxis}
\end{tikzpicture}}
\end{subfloat}
\caption{The $L^2$ errors of various $m$ order modified Hermite-leapfrog methods when applied to the one-dimensional advection equation with unit wave speed. The CFL constant is set to be $C_{CFL}=0.9$ across all orders.}
\label{fig:hermiteAdvecPP}. 
\end{figure}
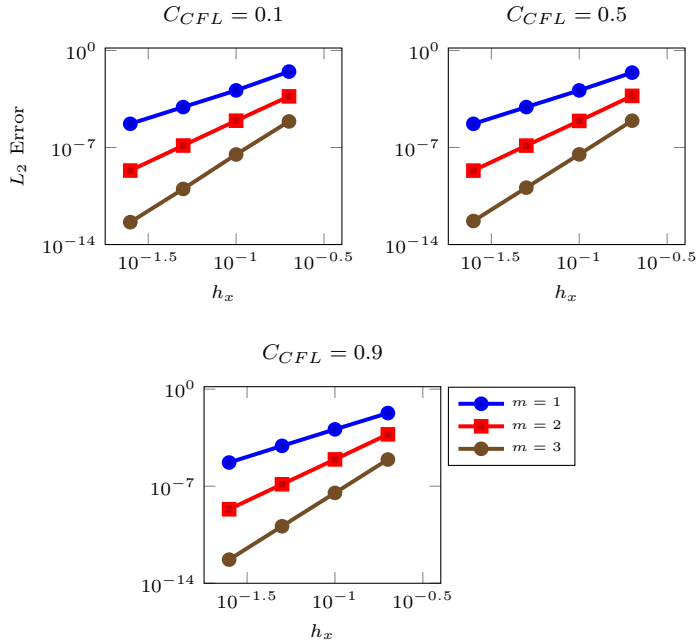

To assess the accuracy and performance of the modified Hermite-leapfrog scheme, numerical experiments are carried out using the one-dimensional advection equation. For these experiments the domain is chosen to be the bi-unit interval and the analytic solution is chosen to be
\[
u(x,t) = \sin(3 \pi(x-t)). 
\]
Figure \ref{fig:hermiteAdvecPP} reports the observed accuracy for orders $m=1,2,3$ while Table \ref{table:advectionV2Rates} reports the observed rates of convergence, notably we observe consistent $O(h^{2m+2})$ rates of convergence.

\begin{table}[h!]
\centering 
\scalebox{0.9}{
\begin{tabular}{|c||c|c|c|c|c|c|c|c|c|}
\hline
& \multicolumn{3}{c|}{$C_{CFL}=0.1$} &  \multicolumn{3}{c|}{$C_{CFL}=0.5$}  &  \multicolumn{3}{c|}{$C_{CFL}=0.9$} \\ 
\hline
Order - $m$ & $1$ & $2$ & $3$ & $1$ & $2$ & $3$  & $1$ & $2$ & $3$ \\ 
\hhline{|=|=|=|=|=|=|=|=|=|=|}
Modified Hermite-leapfrog  & 4.15 & 5.92 & 8.06 & 4.08 & 5.96 & 8.00 & 3.94 & 5.98 & 7.99\\    
\hline
\end{tabular}                                               
}
\caption{Observed $L^2$ rates of convergence when using the modified Hermite-leapfrog scheme to solve the advection equation with unit wave speed.}
\label{table:advectionV2Rates}
\end{table}

\subsubsection{Modified Hermite-leapfrog: Pressure-Velocity Equations}
\begin{figure}[h!]
\centering
\begin{subfloat}{
\begin{tikzpicture}
\begin{loglogaxis}[	
    legend style={font=\tiny},
	scale=0.37,
    title={$m=1$},
    xlabel={$h_x$},
    ylabel={$L_{2}$ Error},
    xmin=0.02,xmax=0.24,
    ymin=0.000001,ymax=3,
    every axis plot/.append style={line width= 1.5pt},
]
\addplot table {HPPComparison_T4_13/HPPRateN1_C9.dat};
\addplot table {HPPComparison_T4_13/LFRateN1_C9.dat};
\addplot table {HPPComparison_T4_13/dualRateN1_C9.dat};
\end{loglogaxis}
\end{tikzpicture}}
\end{subfloat}
\begin{subfloat}{
\begin{tikzpicture}
\begin{loglogaxis}[	
    legend style={font=\tiny},
	scale=0.37,
    title={$m=2$},
    xlabel={$h_x$},
    xmin=0.02,xmax=0.24,
    ymin=0.00000000000001,ymax=3,
    every axis plot/.append style={line width= 1.5pt},
]
\addplot table {HPPComparison_T4_13/HPPRateN2_C9.dat};
\addplot table {HPPComparison_T4_13/LFRateN2_C9.dat};
\addplot table {HPPComparison_T4_13/dualRateN2_C9.dat};
\end{loglogaxis}
\end{tikzpicture}}
\end{subfloat}
\begin{subfloat}{
\begin{tikzpicture}
\begin{loglogaxis}[	
    legend style={font=\tiny},
	scale=0.37,
    title={$m=3$},
    xlabel={$h_x$},
    ylabel={$L_{2}$ Error},
    every axis plot/.append style={line width= 1.5pt},
    xmin=0.02,xmax=0.24,
    ymin=0.00000000000001,ymax=1e-2,
    legend pos= outer north east,
    legend entries={Modified Hermite-leapfrog,Hermite-leapfrog,Dual Hermite},
]
\addplot table {HPPComparison_T4_13/HPPRateN3_C9.dat};
\addplot table {HPPComparison_T4_13/LFRateN3_C9.dat};
\addplot table {HPPComparison_T4_13/dualRateN3_C9.dat};
\end{loglogaxis}
\end{tikzpicture}}
\end{subfloat}
\caption{ $L^2$ errors for the modified Hermite-leapfrog, Hermite-leapfrog, and Dual Hermite scheme when applied to one-dimensional pressure velocity system.}
\label{fig:wavePP}. 
\end{figure}
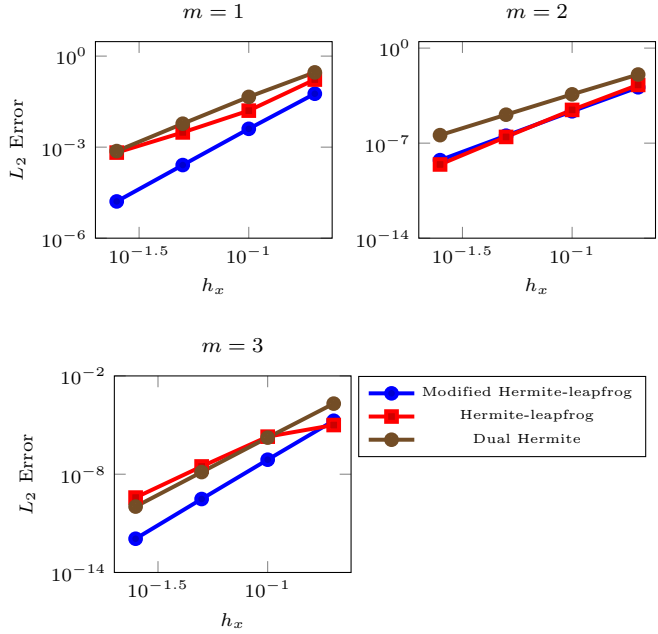

We continue our numerical experiments by applying the modified Hermite-leapfrog method to the pressure-velocity system. The drawback of this approach is that it requires to discretize the pressure and velocity on both the primary and dual grid. Approximations on both grids are necessary as the Cauchy-Kovalesky recurrence relation exchanges even time derivatives of the pressure term for spatial derivatives of the pressure term (odd time derivatives are exchanged for spatial derivatives of the velocity term). For these numerical experiments the analytic solution for the pressure term is chosen to be
\[
p(x,t) = \cos(3\pi t) \sin(3\pi x).
\] 
The system is assumed to have unit wave speed and the solution is propagated to a final time of $T=4.13$. Figure \ref{fig:wavePP} reports the observed accuracy of the modified Hermite-leapfrog scheme under the $L^2$ norm. Table \ref{fig:wavePPRates} reports the observed rates of convergence. The rates of convergence remain consistent with results observed with the advection equation, namely $O(h^{2m+2})$. Notably, the modified Hermite-leapfrog scheme is the most accurate scheme. 

\begin{table}[h!]
\centering 
\scalebox{0.8}{
\begin{tabular}{|c||c|c|c|c|c|c|c|c|c|}
\hline
& \multicolumn{3}{c|}{$C_{CFL}=0.1$} &  \multicolumn{3}{c|}{$C_{CFL}=0.5$}  &  \multicolumn{3}{c|}{$C_{CFL}=0.9$} \\ 
\hline
Order - $m$ & $1$ & $2$ & $3$ & $1$ & $2$ & $3$  & $1$ & $2$ & $3$ \\ 
\hhline{|=|=|=|=|=|=|=|=|=|=|}
Modified Hermite-leapfrog & 4.03 & 6.07 & 7.99 & 3.98 &  5.98 & 8.04 & 3.93 &  5.93 & 7.97  \\ \hline
Hermite-leapfrog & 1.04 & 5.95 & 5.07 & 1.07 & 5.95 & 5.64 & 1.21 & 5.98 & 6.50 \\ \hline
Dual-Hermite & 1.89 & 4.95 & 7.08 & 2.56 & 4.99 & 6.93 & 2.87 & 4.94 & 6.96 \\ 
\hline
\end{tabular} 
}
    \caption{Comparison of the $L^2$ rates of convergence for the modified Hermite-leapfrog, Hermite-leapfrog, and the Dual Hermite method when applied to the one-dimensional pressure-velocity system for orders $m=1,2,3$.}
    \label{fig:wavePPRates}
\end{table} 

\section{Extending the method to higher dimensions}
The Hermite-leapfrog method is easily extended to higher dimensions by means of a tensor product construction. For completeness we describe a two-dimensional extension; a three dimensional version follows a similar construction. To aid in the description of the method we consider the pressure-velocity wave system in two-dimensions 
\begin{align}
    \frac{\partial p}{\partial t} = - \frac{\partial v}{\partial x} - \frac{\partial u}{\partial y}, \\
    \frac{\partial v}{\partial t} = - \frac{\partial p}{\partial x}, \quad \frac{\partial u}{\partial t} = - \frac{\partial p}{\partial y}.
\end{align}
For simplicity we assume the initial conditions are defined at a half-step apart
\begin{align}
p(\mathbf{x},t) = 0, \quad v(\mathbf{x},{\Delta t}/{2}) = g (\mathbf{x}), \quad u (\mathbf{x},{\Delta t}/{2}) = h(\mathbf{x}). 
\end{align}
As in the one-dimensional version, the discretization of the equations occurs over two grids. The primary grid, $\Omega$, serves to hold approximations of the pressure variable. The second dual grid $\tilde{\Omega}$ is introduced in order to store approximations of the velocity fields.  The construction of each grid is performed by taking the tensor product of one-dimensional grids.
Assuming a grid spacing of $h_x$ and $h_y$ in the $x$ and $y$ directions respectively, the solution at each node $(x_j, y_j)$ is represented by the following tensor polynomial 
\[
u(\mathbf{x}) \approx u_{i,j}(\mathbf{x}) = \sum^{m}_{i=0} \sum^{m}_{j=0} \tilde{u}_{i,j} \left( \frac{x-x_i}{h_x} \right)^i \left(\frac{y-y_i}{h_y} \right)^j. 
\]
We emphasize that Hermite-leapfrog methods maintain approximations of the pressure solely on the primary grid. Approximations of the velocity variables are maintained on the dual grid; this is notably different than traditional staggered time domain finite difference methods which discretize variables on different grids. The time-stepping algorithm for the two-dimensional Hermite-leapfrog method is derived in a similar manner as in the one-dimensional case. Here we consider two-dimensional analogues of Equations \ref{eq:option1} and \ref{eq:option2}

\begin{align}
\frac{\partial^{s+l} p(x,y,t+{\Delta t}/{2})}{\partial x^s \partial y^l } &= \sum_{r=0} \frac{1}{r!} \left({\Delta t}/{2}\right)^{r} \frac{\partial^{r+s+l} p(x,y,t)}{\partial t^r \partial x^s \partial y^l } \label{eq:option2D1}, \\
\frac{\partial^{s+l} p(x,y,t-{\Delta t}/{2})}{\partial x^s \partial y^l} &= \sum_{r=0} \frac{1}{r!} \left({\Delta t}/{2}\right)^{r} \frac{\partial^{r+s+l} p(x,y,t)}{\partial t^r \partial x^s \partial y^l }. \label{eq:option2D2} 
\end{align}
As before the update for the pressure term is derived by subtracting Equations \ref{eq:option2D1} and \ref{eq:option2D2} 
\begin{align}
\frac{\partial^{s+l} p(x,y,t+{\Delta t}/{2})}{\partial x^s \partial y^l } - \frac{\partial^{s+l} p(x,y,t-{\Delta t}/{2})}{\partial x^s \partial y^l } = \label{eq:lf1update2D} \\ -\sum_{r=1,odd} \frac{1}{r!} \left({\Delta t}/{2}\right)^{r} \frac{\partial^{r+s+l} v_x(x,y,t)}{\partial x^{s+r} \partial y^{l+r}}  \nonumber \\ 
-\sum_{r=1,odd} \frac{1}{r!} \left({\Delta t}/{2}\right)^{r} \frac{\partial^{r+s+l} v_y(x,y,t)}{\partial x^{s+r} \partial y^{l+r}}  \nonumber.
\end{align}
where the indices $s,l$ correspond to the order of spatial derivatives, i.e. $s,l \in \{0,\dots, m\}$, the formulas for the velocity fields are derived analogously. As in the one-dimensional case; spatial derivatives are approximated by means of Hermite-Birkhoff interpolation. Given the tensor structure of the polynomial, the interpolant enables a dimension by dimension reconstruction. Figure \ref{fig:hermitePlus2D} illustrates the staggering of the nodes as well as the reconstruction procedure; we refer the reader to \cite{vargas2017hermite} for further details the higher dimensional reconstruction procedure. 

\begin{figure}[h!]
\centering
\subfloat[]{\scalebox{0.4}{\includegraphics{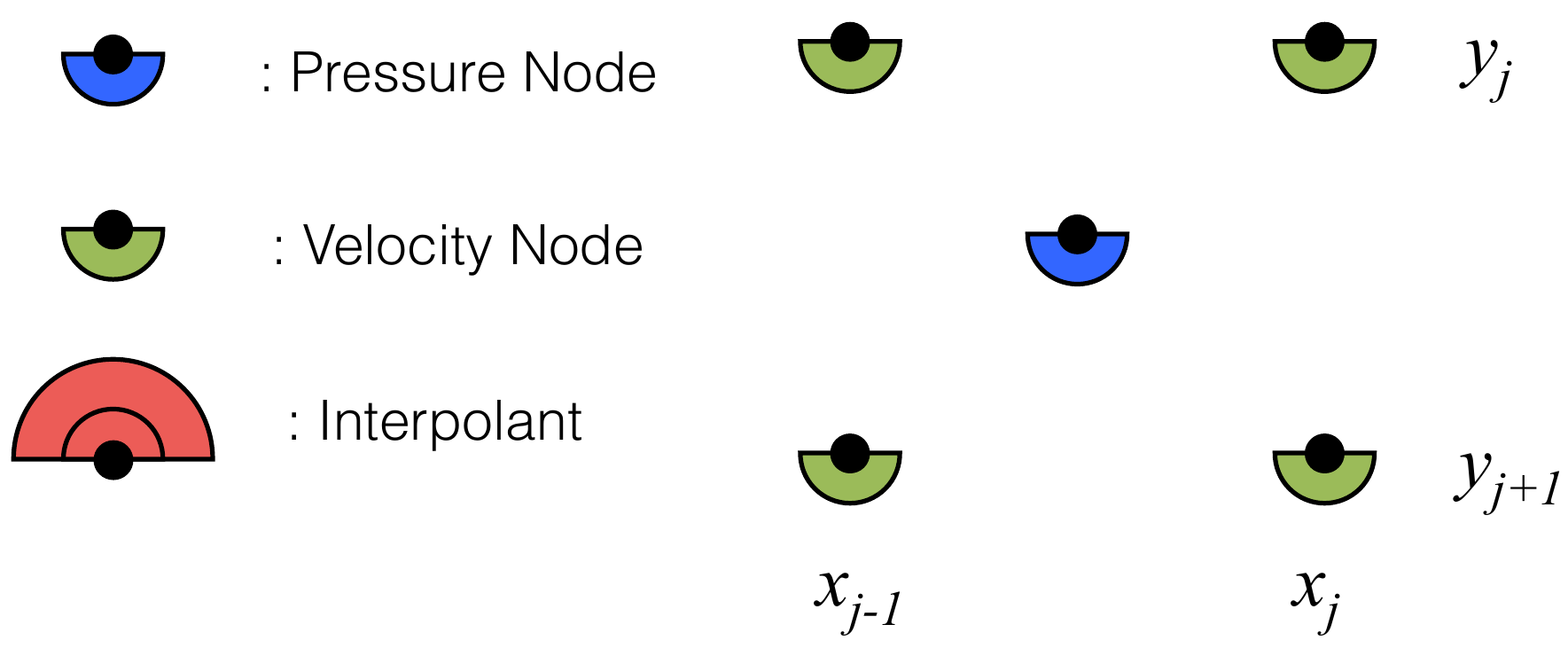}}} \quad
\subfloat[]{\scalebox{0.4}{\includegraphics{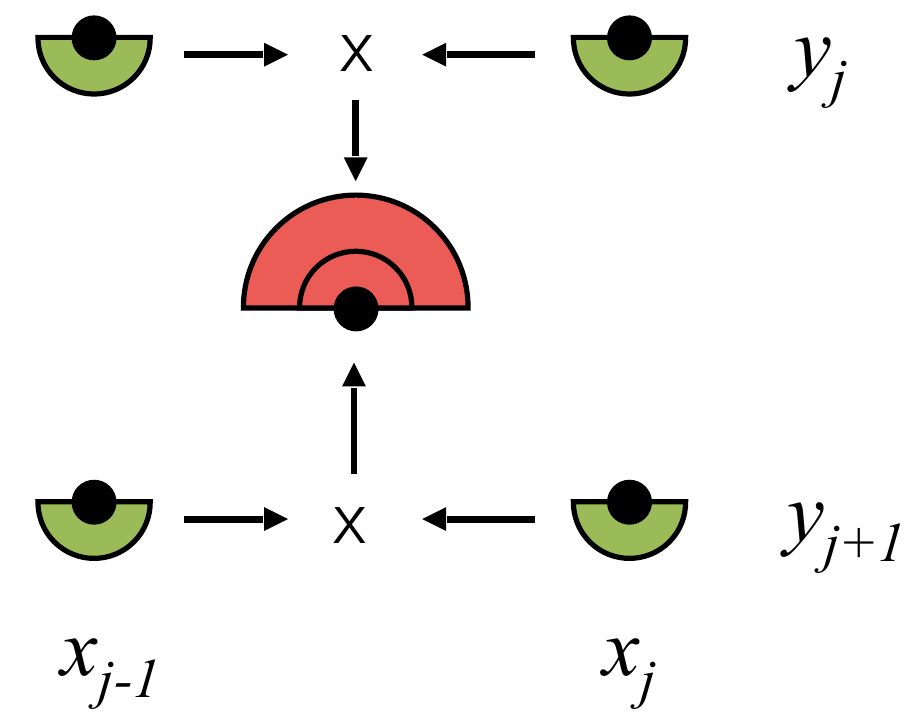}}}
\caption{Figure (a) illustrates an aerial view of the Hermite-leapfrog discretization for the pressure-velocity equations. We emphasize that the nodes are staggered a half-step apart. Figure (b) illustrates a dimension by dimension reconstruction of the velocity interpolant. The interpolants tensor product structure enables a dimension by dimension reconstruction.}
\label{fig:hermitePlus2D}. 
\end{figure}

\subsection{Numerical experiments in two-dimensions}
To demonstrate the method's efficiency, we solve the acoustic wave equations with unit wave speed on the bi-unit square, $[-1,1] \times [-1,1]$, with periodic boundary conditions. The pressure is chosen to have an analytic solution of
\[
p(x,y,t) = \sin(\pi x) \sin(\pi y) \cos(\sqrt{2 \pi} t). 
\]
Figure \ref{fig:hermitePlus2D} reports observed accuracy in the $L^2$ norm for both the Hermite-leapfrog and Dual Hermite method. Table \ref{table:hermitePlus2D} reports observed rates of convergence with varying CFL constants, notably we observe similar convergence behavior as in the one-dimensional case.  

\begin{figure}[h!]
\centering
\begin{subfloat}{
\begin{tikzpicture}
\begin{loglogaxis}[	
    legend style={font=\tiny},
	scale=0.36, 
    title={$m=0$},
    xlabel={$h_x$},
    ylabel={$L_{2}$ Error},
    xmin=0.01,xmax=0.21,
    ymin=0.0001,ymax=1.5,
    every axis plot/.append style={line width= 1.5pt}
]
\addplot table {hPlusRatesV2_T4_13/dualRateN0_C9.dat};
\addplot table {hPlusRatesV2_T4_13/hPlusRateN0_C9.dat};
\end{loglogaxis}
\end{tikzpicture}}
\end{subfloat}
\begin{subfloat}{
\begin{tikzpicture}
\begin{loglogaxis}[	
    legend style={font=\tiny},
	scale=0.36, 
    title={$m=1$},
    xlabel={$h_x$},
    xmin=0.01,xmax=0.21,
    ymin=1e-5,ymax=1e-2,
    legend pos= outer north east,
    every axis plot/.append style={line width= 1.5pt},
    legend entries={Dual Hermite, Hermite-leapfrog},
]
\addplot table {hPlusRatesV2_T4_13/dualRateN1_C9.dat};
\addplot table {hPlusRatesV2_T4_13/hPlusRateN1_C9.dat};
\end{loglogaxis}
\end{tikzpicture}}
\end{subfloat}
\begin{subfloat}{
\begin{tikzpicture}
\begin{loglogaxis}[	
    legend style={font=\tiny},
	scale=0.36, 
    title={$m=2$},
    xlabel={$h_x$},
    ylabel={$L_{2}$ Error},
    xmin=0.01,xmax=0.21,
    ymin=1e-12,ymax=1.5e-3,
    every axis plot/.append style={line width= 1.5pt}
]
\addplot table {hPlusRatesV2_T4_13/dualRateN2_C9.dat};
\addplot table {hPlusRatesV2_T4_13/hPlusRateN2_C9.dat};
\end{loglogaxis}
\end{tikzpicture}}
\end{subfloat}
\begin{subfloat}{
\begin{tikzpicture}
\begin{loglogaxis}[	
    legend style={font=\tiny},
	scale=0.36, 
    title={$m=3$},
    xlabel={$h_x$},
    xmin=0.01,xmax=0.21,
    ymin=1e-15, ymax=1e-3,
    legend pos= outer north east,
    every axis plot/.append style={line width= 1.5pt},
    legend entries={Dual Hermite, Hermite-leapfrog},
]
\addplot table {hPlusRatesV2_T4_13/dualRateN3_C9.dat};
\addplot table {hPlusRatesV2_T4_13/hPlusRateN3_C9.dat};
\end{loglogaxis}
\end{tikzpicture}}
\end{subfloat}
\caption{ $L^2$ errors for the modified Hermite-leapfrog, Hermite-leapfrog, and Dual Hermite scheme when applied to two-dimensional pressure velocity system.}
\label{fig:hermitePlus2D}. 
\end{figure}
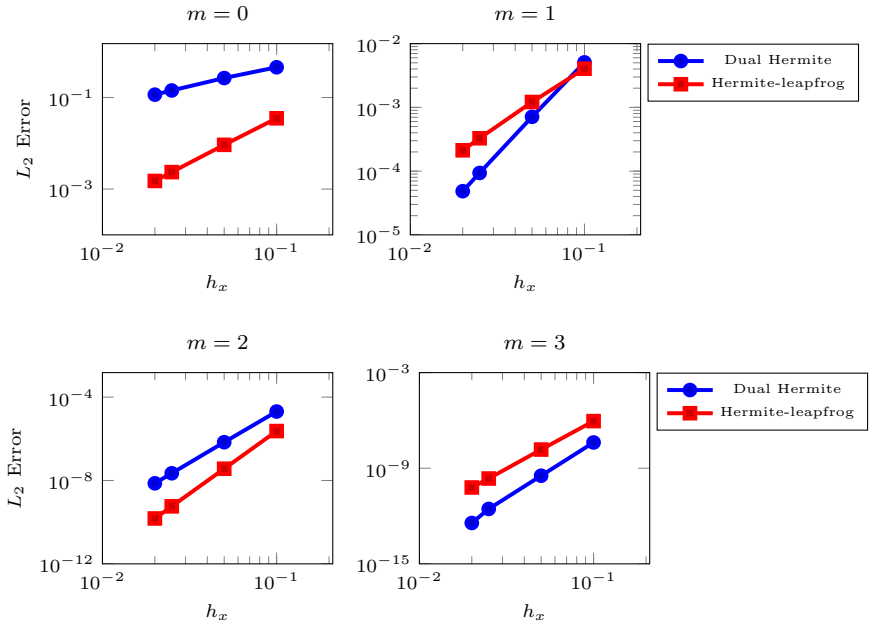

\begin{table}[h!]
\centering
\scalebox{0.8}{   
\begin{tabular}{|l|l|l|l|l|l|l|l|l|l|l|l|l|}
\hline
             & \multicolumn{4}{l|}{$C_{CFL} = 0.1$} & \multicolumn{4}{l|}{$C_{CFL} = 0.5$} & \multicolumn{4}{l|}{$C_{CFL} = 0.9$} \\ \hline
Order - $m$            & 0     & 1    & 2    & 3    & 0     & 1     & 2    & 3    & 0     & 1    & 2    & 3    \\ \hline
Hermite-leapfrog & 1.74 & 1.76 & 5.68 & 6.84 & 1.75 & 1.76 & 5.71 & 7.45 & 1.86 & 1.88 & 6.01 & 6.74 \\ \hline
Dual-Hermite & -  & 2.37 & 4.76 & 6.83 & -  & 2.44 & 4.78 & 6.76 & - & 2.67 & 4.80 & 6.77 \\ \hline
\end{tabular}
}
\caption{Comparison of $L^2$ rates of convergence for the $m$ order Hermite-leapfrog, and the Dual Hermite method when applied to the one-dimensional pressure-velocity system.}
\label{table:hermitePlus2D}
\end{table}

\subsection{Imposing reflective boundary conditions}
\label{sec:reflec}
\begin{figure}[h!]
\centering
\subfloat[]{\scalebox{0.3}{\includegraphics{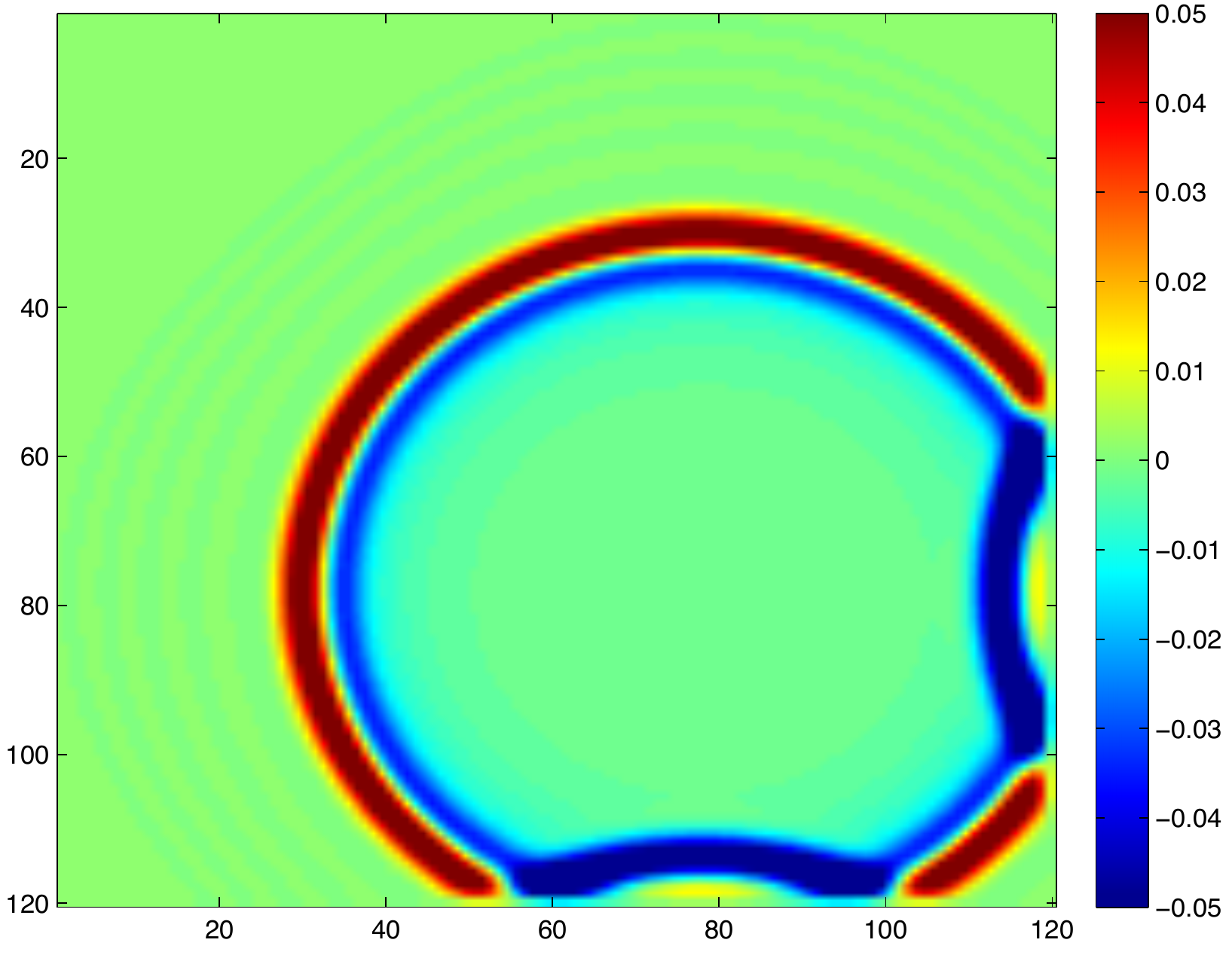}}} \quad
\subfloat[]{\scalebox{0.3}{\includegraphics{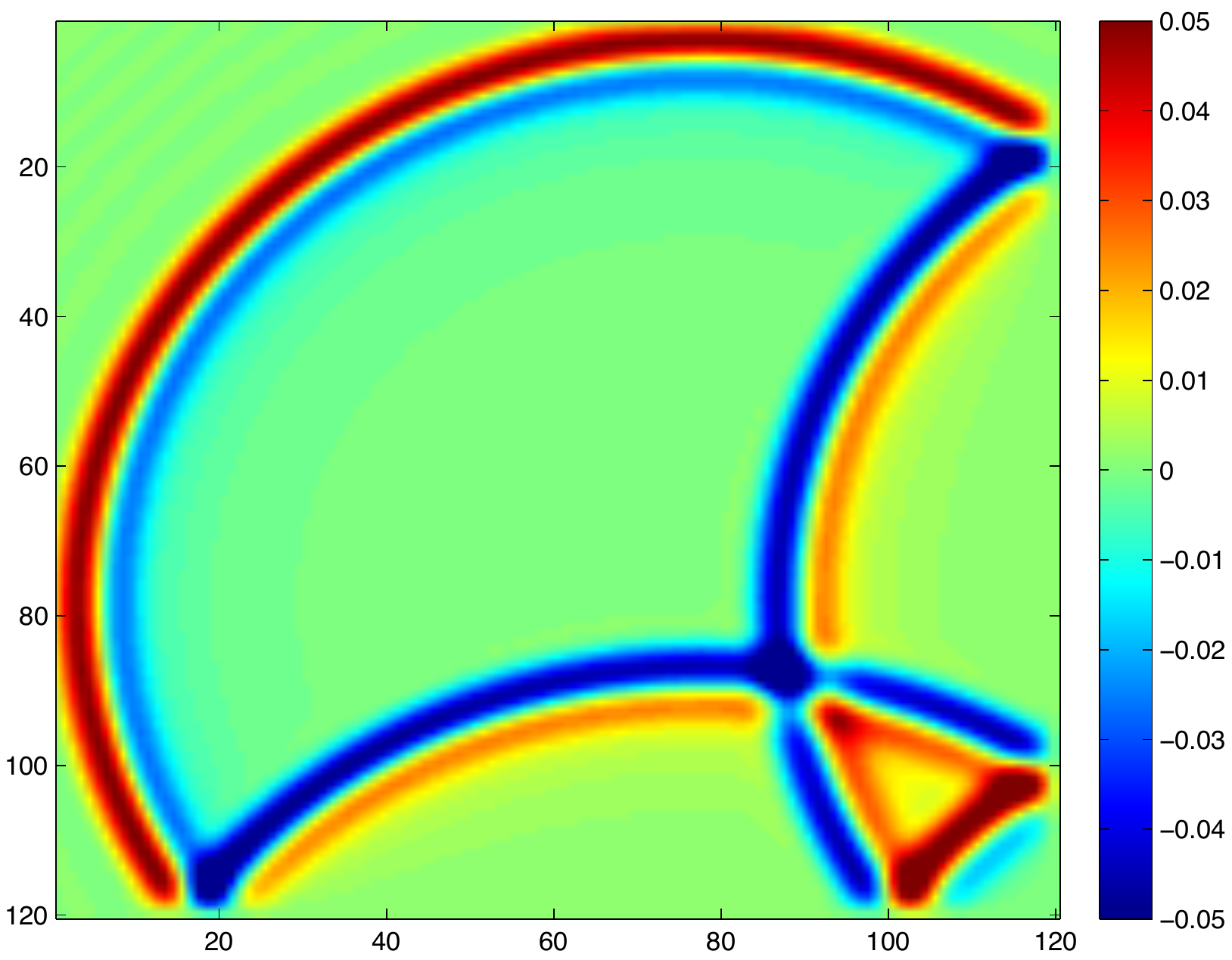}}} 
\caption{Figure (a) and (b) illustrate a wave reflected on the South-East corner as propagated by the Hermite-leapfrog scheme. Imposing reflective boundary conditions is accomplished by adding ghost nodes on the dual grid. The ghost points values are chosen to mirror the polynomials such that they are even or odd around the boundary and thus agreeing with the boundary conditions specified by the boundary conditions (Equations \ref{eq:Boundary}) .  }
\label{fig:zeroWall}. 
\end{figure}

As a second example, we consider zero Dirichlet boundary conditions on the pressure variable. The complete set of boundary conditions may be derived by insisting the solution of the pressure term is zero at the boundary and differentiating in time; yielding
\begin{align}
 p(x,y,t) &= 0, \quad x = -1,1, \quad    y = -1,1 \label{eq:Boundary} \\ 
 u(x,y,t) &= 0, \quad u_{y}(x,y,t) = 0, \quad y = -1, 1 \nonumber \\ 
 v(x,y,t) &= 0, \quad v_{x}(x,y,t) = 0, \quad x = -1, 1 \nonumber.
\end{align}
The boundary conditions are imposed prior to propagating the pressure solution by filling in ghost point values of the velocity variables at the edge of the domain. The ghost point values are chosen to mirror the polynomials such that they are either odd around the boundary and thus agree with the homogeneous Dirichlet boundary conditions or so that they are even and satisfy homogeneous Neumann conditions. To illustrate reflecting boundaries we initialize the pressure variable with the following Gaussian pulse \begin{align}
p(x,y,0) = \exp \left( \frac{(x-0.3)^2 + (y-0.3)^2}{0.002} \right). \nonumber
\end{align}
Figure \ref{fig:zeroWall} illustrates a reflection of the propagating wave on the South-East corner.

\subsection{Numerical experiments with Maxwell's equations}
\begin{figure}
\centering
\begin{tikzpicture}
\begin{loglogaxis}
[	
    legend style={font=\tiny},
	scale=0.6, 
    title={$CFL=0.8$},
    xlabel={$h_x$},
    xmin=2.8986e-02,xmax=1.0e-01,
    ymin=1e-10,ymax=2,
    legend pos= outer north east,
    every axis plot/.append style={line width= 1.5pt},
    legend entries={$m=3$,$m=4$,$m=5$,$m=6$},
]
\addplot table {maxwells2/m_3_08.dat};
\addplot table {maxwells2/m_4_08.dat};
\addplot table {maxwells2/m_5_08.dat};
\addplot table {maxwells2/m_6_08.dat};
\end{loglogaxis}
\end{tikzpicture}
\caption{ $L^2$ errors for the Hermite-leapfrog method applied to Maxwell's equation.}
\label{fig:maxwells}
\end{figure}
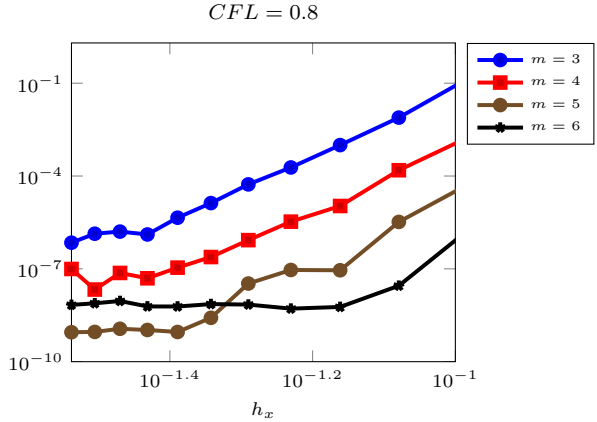

To illustrate our methods ability to resolve high frequency waves,
we carry out numerical experiments with Maxwell's equation in the time domain. In particular we consider the transverse magnetic form with unit material coefficients
\begin{align}
    \frac{\partial H^{x}} {\partial t} &= -\frac{\partial E^{z} } {\partial y},  \\
    \frac{\partial H^{y}} {\partial t} &= \frac{\partial E^{z} } {\partial x},  \nonumber \\
    \frac{\partial E^{z}} {\partial t} &= \frac{\partial H^{y} } {\partial x} - \frac{\partial H^{x}}{\partial y}. \nonumber
\end{align}
Here $(H^{x}, H^{y})$ corresponds to the magnetic fields, and $E^{z}$ is the electric field. We take our domain of interest to be the bi-unit square, $\Omega = [-1, 1] \times [-1, 1]$, and the analytic solution to be
\begin{align}
    H^{x}(x,y,t) &= -\frac{\omega_y}{\omega_t} \sin(\omega_x x) \cos(\omega_y y) \sin(\omega_t t),  \\
    H^{y}(x,y,t) &= \frac{\omega_x}{\omega_t} \cos(\omega_x x) \sin(\omega_y y) \sin(\omega_t t), \nonumber \\
    E^{z}(x,y,t) &= \sin(\omega_x x) \cos(\omega_y y) \sin(\omega_t t). \nonumber 
\end{align}
The parameters are chosen to be
$\omega_x = 8 \pi, \omega_y = 8 \pi$, and $\omega_t = \sqrt{  \omega_x^2 + \omega_y^2  }$.
The boundary conditions are imposed by the mirroring technique discussed in \ref{sec:reflec}. To illustrate the resolving power of a high order numerical method, we choose the polynomial degree of the Hermite-leapfrog method to be $m=3,4, 5$, and $6$. Figure \ref{fig:maxwells} reports the error under the $L_2$ norm. At high orders we observe slightly more variation of rates of convergence and sensitivity to round-off, yet the method is still numerically stable. 

\section{Acceleration on graphics processing units}
Lastly, we demonstrate that the Hermite-Leapfrog method is well suited for the graphics processing unit. In these numerical experiments we solve the three-dimensional acoustic wave equations
\begin{align}
    \frac{\partial p}{\partial t} = - \frac{\partial v}{\partial x} - \frac{\partial u}{\partial y} - \frac{\partial w}{\partial z} \\
    \frac{\partial v}{\partial t} = - \frac{\partial p}{\partial x}, \quad \frac{\partial u}{\partial t} = - \frac{\partial p}{\partial y}, \quad \frac{\partial w}{\partial t} = \frac{\partial p}{\partial z} \nonumber. 
\end{align}
on the unit cube, $[-1,1] \times [-1, 1]$ with periodic boundary conditions. This section serves as an extension of the work carried out by Vargas et at. \cite{vargas2016gpu} wherein Hermite methods were tailored to the graphics processing unit. More precisely we revisit the two approaches outlined in \cite{vargas2016gpu}. The first is a split approach in two kernels are employed. The first kernel reconstructions approximations of high order derivatives via Taylor series while the second kernel then applies leapfrog time-stepping to propagate the solution. The second approach combines the interpolation and evolution procedures into a monolithic kernel.

Our numerical experiments are carried out on a single node of the Southern Methodist University cluster ``ManeFrame II". A single node consist of a NVIDIA Tesla P100-PCIE-16GB graphics card. The hardware has a theoretical bandwidth of 732 GB/s and potentially perform 4700 Gigaflops/sec in double precision. To obtain estimates of kernel bandwidth and gflops we employ the NVIDIA profiler which reports effective arithmetic bandwidth by flop\_count\_dp and memory bandwidth as the sum of \textit{dram\_read\_throughput} and \textit{dram\_write\_throughput}. Here bandwidth corresponds to the sum of bytes read and written to global memory by a GPU kernel. Our CUDA kernels are generated by using the OCCA language \cite{medina2015okl}. Lastly we compare time to solution with a dual socket Xeon CPU E5-2695 (2.10 GHz) as found on a Maneframe II node.

\begin{figure}[h!]
\centering
\subfloat{
\begin{tikzpicture}
\begin{axis}[
	legend style={font=\tiny},
	width=.45\textwidth,
	xmin=0.0,xmax=4,
	xtick={1,2,3},
	ymin=0,
	ymax=750,
	ylabel=Bandwidth (GB/s),
	xlabel=Degree $m$,
	ybar=4pt,
	bar width=5pt,
	xmajorgrids=true,
	ymajorgrids=true,
	grid style=dashed,
	legend pos=north west,
	title=Bandwidth
]

\addplot 
    coordinates {(1,433) (2,438) (3,386)}; 

\addplot 
    coordinates {(1,445) (2,335) (3,344)};
	
\addplot 
    coordinates{(1,380) (2,256) (3,245)};
	



\end{axis}
\end{tikzpicture}
}
\subfloat{
\begin{tikzpicture}
\begin{axis}[
	legend style={font=\tiny},
    xmin=0.0,xmax=4,
	width=.45\textwidth,
	xtick={1,2,3},
	ylabel=GFLOP/s,
	xlabel=Degree $m$,
	ymin=50,
	ymax=1500,
	ybar=4pt,
	bar width=5pt,
	xmajorgrids=true,
	ymajorgrids=true,
	legend pos=north west,
	grid style=dashed,
	title=GFLOP/s,
	legend pos=outer north east
]

\addplot 
    coordinates {(1,775) (2,1222) (3,1334)};

\addplot 
    coordinates{(1,228) (2,284) (3,188)};
	
\addplot 
   coordinates{(1,1226) (2,1154) (3,596.9)};

\legend{Inter, Evo Pre, Evo Vel}
\end{axis}
\end{tikzpicture}
}

\caption{Here we report the Bandwidth and GFLOP of a two kernel approach as outlined by Vargas and co-authors in \cite{vargas2016gpu}. }\label{fig:HermiteYeeGPU1}
\end{figure}
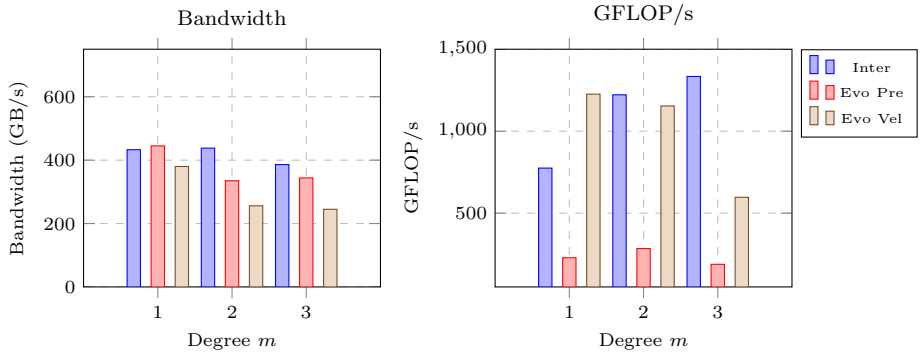

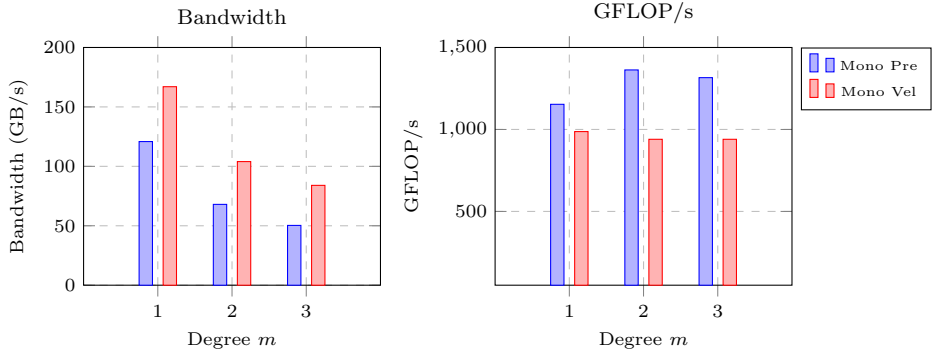
\begin{figure}[h!]
\centering
\subfloat{
\begin{tikzpicture}
\begin{axis}[
	legend style={font=\tiny},
	width=.45\textwidth,
	xmin=0.0,xmax=4,
	xtick={1,2,3},
	ymin=0,
	ymax=200,
	ylabel=Bandwidth (GB/s),
	xlabel=Degree $m$,
	ybar=4pt,
	bar width=5pt,
	xmajorgrids=true,
	ymajorgrids=true,
	grid style=dashed,
	legend pos=north west,
	title=Bandwidth
]

\addplot 
    coordinates{(1,120.84) (2,68) (3,50.3)};

\addplot 
    coordinates{(1,167) (2,104) (3,84)};



\end{axis}
\end{tikzpicture}
}
\subfloat{
\begin{tikzpicture}
\begin{axis}[
	legend style={font=\tiny},
    xmin=0.0,xmax=4,
	width=.45\textwidth,
	xtick={1,2,3},
	ylabel=GFLOP/s,
	xlabel=Degree $m$,
	ymin=50,
	ymax=1500,
	ybar=4pt,
	bar width=5pt,
	xmajorgrids=true,
	ymajorgrids=true,
	legend pos=north west,
	grid style=dashed,
	title=GFLOP/s,
	legend pos=outer north east
]

\addplot 
    coordinates{(1,1153.38) (2,1363) (3,1316)};

\addplot 
    coordinates{(1,987) (2,940) (3,940)};

\legend{Mono Pre, Mono Vel}
\end{axis}
\end{tikzpicture}
}

\caption{Here we report the Bandwidth and GFLOP of using a monolithic kernel approach as outline by Vargas and co-authors in \cite{vargas2016gpu}.}\label{fig:HermiteYeeGPU2}
\end{figure}

In studying GFLOP and effective bandwidth it is important to note that the devices theoretical capabilities are typically difficult to achieve in practice. To estimate a more ``realistic'' bandwidth we consider a simple vector copy wherein the entries of vector $x$ are copied onto vector $y$. For sufficiently large $N$ a streaming bandwidth of 534 GB/s was observed on the graphics card. Although this is not indicative of the peak performance it does serve as a representative of achievable peak performance numbers. Tables \ref{fig:HermiteYeeGPU1} and Tables \ref{fig:HermiteYeeGPU2} report the observed GFLOP and bandwidth performance. 

In our numerical experiments it is clear that the interpolation procedure dominates the floating point operations. This is not too surprising as the operations are carried out as a series of matrix-matrix multiplications (tensor contractions). Similarly the kernel which evolves the velocity variables, $V_x, V_y, V_z$, has a high number of floating point operations as it carried out the time-stepping for three variables.

As OCCA enables code generation on various platforms we conclude our accelerator studies by presenting Table \ref{table:cpuvgpuRuntime} which compares time per iteration across a P100 and the dual 16-core Xeon CPUs. In these experiments an iteration refers to propagating the pressure and velocity fields.

\begin{table}[h!]
\centering
\label{table:timeToSolution}
\begin{tabular}{|l|c|c|c|c|}
\hline
Order of the method - $m$        &  $m=1$  & $m=2$ & $m=3$ \\ \hline
OCCA::OpenMP  &  1.08 sec &  1.30 sec &  1.74 sec \\ \hline
OCCA::CUDA - Single Kernel   & 0.060 sec &  0.07 sec &  0.10 sec \\ \hline
OCCA::CUDA - Two Kernels    & 0.098 sec & 0.19 sec &  0.15 sec\\ \hline
\end{tabular}
\caption{Comparison of time per iteration of Hermite kernels executed on the GPU and CPU. For orders $m=1,2,3$ the number of grid points were chosen to be 200, 130, 100 points per Cartesian direction. Noticeably a single GPU kernel offers a better time to solution with advantage of reliving the need to store the interpolate on account of less data movement. These results are consistent with what has been observed in \cite{vargas2016gpu}.}
\label{table:cpuvgpuRuntime}
\end{table}

\section{Summary}
We have presented a variation of the classic Hermite method which uses leap-frog time-stepping to advance the solution. The new Hermite-leapfrog method is demonstrated to be numerically stable, high-order in both time and space, and may time-step the solution independent of order. A detailed description of the method is provided in one and two dimensions, as well as a techniques for incorporating spatially varying coefficients and reflective boundary conditions. Numerical experiments suggest that the rate of convergence is dependent on whether the method is of odd or even order. To address the variation we introduce a modification that achieves consistent $O(h^{2m+2})$ rates of convergence. 

Lastly, we accelerate the method using a graphics processing unit. We find that a monolithic kernel is ideal as it eliminates the need to explicitly store the interpolant and can provide a comparable time to solution to using two kernels which achieve a higher bandwidth. As a future research direction we plan to provide a theoretical justification for the exceptional convergence rates of the modified Hermite-leapfrog method.

\begin{acknowledgements}
The authors would like to thank Daniel Appel\"{o} for the fruitful conversations. TH is supported in part by  NSF Grant DMS-1418871. Any conclusions or recommendations expressed in this paper are those of the author 
and do not necessarily reflect the views of the NSF.J C is supported by NSF grants DMS-1719818 and DMS-1712639. This work was performed under the auspices of the U.S. Department of Energy by Lawrence Livermore National Laboratory under Contract DE-AC52-07NA27344. LLNL-JRNL-757049.
\end{acknowledgements}

\bibliographystyle{spmpsci}      
\bibliography{main.bbl}   

\begin{thebibliography}{10}
\providecommand{\url}[1]{{#1}}
\providecommand{\urlprefix}{URL }
\expandafter\ifx\csname urlstyle\endcsname\relax
  \providecommand{\doi}[1]{DOI~\discretionary{}{}{}#1}\else
  \providecommand{\doi}{DOI~\discretionary{}{}{}\begingroup
  \urlstyle{rm}\Url}\fi

\bibitem{appelo2018hermite}
Appelo, D., Hagstrom, T., Vargas, A.: Hermite methods for the scalar wave
  equation.
\newblock arXiv preprint arXiv:1802.05246  (2018)

\bibitem{appelo2011hermite}
Appel{\"o}, D., Inkman, M., Hagstrom, T., Colonius, T.: {Hermite methods for
  aeroacoustics: Recent progress}.
\newblock In: 17th AIAA/CEAS Aeroacoustics Conference (32nd AIAA Aeroacoustics
  Conference), Portland, Oregon (2011)

\bibitem{bencomo2015discontinuous}
Bencomo, M.J.: {Discontinuous Galerkin and finite difference methods for the
  acoustic equations with smooth coefficients}.
\newblock Ph.D. thesis, Rice University (2015)

\bibitem{chen2012p}
Chen, R., Hagstrom, T.: {P-adaptive Hermite methods for initial value
  problems}.
\newblock ESAIM: Mathematical Modelling and Numerical Analysis \textbf{46}(3),
  545--557 (2012)

\bibitem{chen2014hybrid}
Chen, X.R., Appel{\"o}, D., Hagstrom, T.: {A hybrid Hermite--discontinuous
  Galerkin method for hyperbolic systems with application to Maxwellʼs
  equations}.
\newblock Journal of Computational Physics \textbf{257}, 501--520 (2014)

\bibitem{gauthier1986two}
Gauthier, O., Virieux, J., Tarantola, A.: {Two-dimensional nonlinear inversion
  of seismic waveforms: Numerical results}.
\newblock Geophysics \textbf{51}(7), 1387--1403 (1986)

\bibitem{goodrich2006hermite}
Goodrich, J., Hagstrom, T., Lorenz, J.: {Hermite methods for hyperbolic
  initial-boundary value problems}.
\newblock Mathematics of computation \textbf{75}(254), 595--630 (2006)

\bibitem{kornelus2018flux}
Kornelus, A., Appel{\"o}, D.: Flux-conservative hermite methods for simulation
  of nonlinear conservation laws.
\newblock Journal of Scientific Computing \textbf{76}(1), 24--47 (2018)

\bibitem{kowalczyk2010comparison}
Kowalczyk, K., Van~Walstijn, M.: {A comparison of nonstaggered compact FDTD
  schemes for the 3D wave equation}.
\newblock In: Acoustics Speech and Signal Processing (ICASSP), 2010 IEEE
  International Conference on, pp. 197--200. IEEE (2010)

\bibitem{levander1988fourth}
Levander, A.R.: {Fourth-order finite-difference P-SV seismograms}.
\newblock Geophysics \textbf{53}(11), 1425--1436 (1988)

\bibitem{leveque2007finite}
LeVeque, R.J.: {Finite difference methods for ordinary and partial differential
  equations: steady-state and time-dependent problems}.
\newblock SIAM (2007)

\bibitem{medina2015okl}
Medina, D.: {OKL: A unified language for parallel architectures}.
\newblock Ph.D. thesis, Rice University (2015)

\bibitem{vargas2017hermite}
Vargas, A.: {Hermite methods for the simulation of wave propagation}.
\newblock Ph.D. thesis (2017)

\bibitem{vargas2016gpu}
Vargas, A., Chan, J., Hagstrom, T., Warburton, T.: Gpu acceleration of hermite
  methods for the simulation of wave propagation.
\newblock In: Spectral and High Order Methods for Partial Differential
  Equations ICOSAHOM 2016, pp. 357--368. Springer (2017)

\bibitem{vargas2017variations}
Vargas, A., Chan, J., Hagstrom, T., Warburton, T.: {Variations on Hermite
  methods for wave propagation}.
\newblock Communications in Computational Physics \textbf{22}(2), 303--337
  (2017)

\bibitem{xie2002explicit}
Xie, Z., Chan, C.H., Zhang, B.: An explicit fourth-order staggered
  finite-difference time-domain method for maxwell's equations.
\newblock Journal of Computational and Applied Mathematics \textbf{147}(1),
  75--98 (2002)

\bibitem{yee1966numerical}
Yee, K.: {Numerical solution of initial boundary value problems involving
  Maxwell's equations in isotropic media}.
\newblock IEEE Transactions on antennas and propagation \textbf{14}(3),
  302--307 (1966)

\bibitem{yefet2001staggered}
Yefet, A., Petropoulos, P.G.: A staggered fourth-order accurate explicit finite
  difference scheme for the time-domain maxwell's equations.
\newblock Journal of Computational Physics \textbf{168}(2), 286--315 (2001)

\end{thebibliography}

%
%

\end{document}